\documentclass{LMCS}

\def\doi{9(3:5)2013}
\lmcsheading%
{\doi}
{1--24}
{}
{}
{May\phantom.~19, 2011}
{Aug.~20, 2013}
{}

\usepackage{enumerate}
\usepackage{hyperref}
\usepackage{path}

\usepackage{amssymb}
\usepackage{amsmath}
\usepackage{latexsym}

\newcommand{\hide}[1]{}

\def \dom{{\rm dom}}

\def \In{\mathrel{\subseteq}}
\def \om{{\Sigma^\omega }}
\def \pf{:\hspace{0.6ex}\subseteq \hspace{-0.4ex}}

\def \prf{\sqsubseteq}
\def \range{{\rm range}}
\def \s{{\Sigma^*}}

\def \dm{{\rm dm}}

\def \pref {\sqsubseteq}

\def\IN{{\mathbb{N}}}
\def\IN{{\mathbb{N}}}
\def\IQ{{\mathbb{Q}}}
\def\IR{{\mathbb{R}}}
\def\IR{{\mathbb{R}}}
\def\IZ{{\mathbb{Z}}}
\newcommand{\an}{\ \wedge\ }

\newcommand{\mto}{\rightrightarrows}

\newcommand{\nufs}{{       \nu^{\rm fs}}}

\newcommand{\nui} {{\bigcap\nu^{\rm fs}}}

\newcommand{\pproof}{\noindent{\bf Proof:} }
\newcommand{\qq}{\phantom{.}\hfill$\Box$}

\newcommand{\bb}{ \hspace{-0.76ex}- \hspace{-0.80ex}}% für korrekte Silbentrennung

\def\Ty{T\hspace{-.4ex}y}

\newcommand{\msim}{{/\hspace{-0.7ex}\sim}}

\newcommand{\mmto}{\mbox{
\setlength{\unitlength}{1em}
\begin{picture}(0.4,0)
\makebox(0,0.6){$\mbox{\scriptsize \raisebox{0.083em}{$|$}}   \hspace*{-1.1ex}\mto$}
\end{picture}
}}

\begin{document}

\title[Computably regular spaces]{ Computably regular topological spaces}

\author[K.~Weihrauch]{Klaus Weihrauch} %optional
\address{University of Hagen, Hagen, Germany}    %optional
\email{Klaus.Weihrauch@FernUni-Hagen.de}  %optional
%\thanks{thanks 2, optional.}    %optional

\keywords{computable analysis, computable topology, axioms of separation}
\subjclass{F.0, F.m, G.0, G.m}
\ACMCCS{[{\bf Mathematics of computing}]: Continuous
  mathematics---topology---point set topology}

\begin{abstract}
This article continues the study of computable elementary topology
started by the author and T.~Grubba in 2009 and extends the author's
2010 study of axioms of computable separation. Several computable
$T_3$- and Tychonoff separation axioms are introduced and their
logical relation is investigated. A number of implications between
these axioms are proved and several implications are excluded by
counter examples, however, many questions have not yet been
answered. Known results on computable metrization of $T_3$-spaces from
M.~Schr\"oder (1998) and T.~Grubba, M.~Schr\"oder and the author
(2007) are proved under uniform assumptions and with partly simpler
proofs, in particular, the theorem that every computably regular
computable topological space with non-empty base elements can be
embedded into a computable metric space.  Most of the computable
separation axioms remain true for finite products of spaces.

\end{abstract}
%%%%%%%%%%%%%%%%%%%%%%%%%5

\maketitle

\section{Introduction}\label{seca}

This article continues with the study of computable topology started in \cite{WG09}.
For computable topological spaces (as defined in \cite{WG09}) in \cite{Wei10} we have introduced a number of computable versions of the
topological $T_0$-, $T_1$- and $T_2$-axioms and studied  their relationship.
In this article we define various computable versions of the topological $T_3$-, Tychonoff- and $T_4$-axioms and compare them. Furthermore, we study computable metrization.
For classical topology see, for example, \cite{Eng89}.
In addition to new material we include earlier results from \cite{Sch98,GSW07,GWX07a,XG09} and \cite{Wei09} (in \cite{BHK09}) some of which have been proved under slightly differing assumptions and give some simpler  proofs.

We will use the representation approach of computable analysis \cite{KW85,Wei00,BHW08}. As the basic computability structure we start with computable topological spaces as introduced in \cite{WG09}. Notice that there are other slightly differing not equivalent definitions of  ``computable topological space'' in other publications, in particular in \cite{Wei00}.
We will use the notations and results from \cite{WG09} some of which are mentioned
very shortly in Section~\ref{secb}.

In Section \ref{secf} we introduce axioms for computable $T_2$ (2 axioms, which are  alredy studied in \cite{Wei10}), for computable
$T_3$ (3 axioms), for computable Tychonoff (3 axioms) and for computable $T_4$ and computable Urysohn. We give some examples and prove that the axioms do not depend
on the details of the computable topological space but only on the computability concept defined by it.

In Section~\ref{seco} we prove a number of implications between the introduced axioms.

In Section~\ref{secp} we show by counterexamples that some implications are false. We  summarize the results and list some open problems concerning the implications between the axioms. We also prove that computable $T_3$ and computable Tychonoff as well as their strong versions  are equivalent for computable topological spaces with non-empty base elements.

In Section~\ref{secs} we resume results on computable metrization from \cite{Sch98,GSW07,GWX07a} and prove them under common weak assumptions. In particular we give a considerably simpler proof of the main theorem from \cite{GSW07} on the embedding computable $T_3$-spaces in computable metric spaces.

Each of the introduced computable separation classes is closed under the subspace operations, and most of them are closed under Cartesian product (Section~\ref{secr}).

\section{Preliminaries}\label{secb}

We will use the terminology and abbreviations summarized in \cite[Section~2]{WG09}
and also results from \cite{WG09}.
For further details see \cite{Wei00,Wei08,BHW08}.

Let $\Sigma$ be a finite alphabet such that $0,1\in\Sigma$. By $\s$ we denote the set of
finite words over $\Sigma$ and by $\om$ the set of infinite sequences
$p:\IN\to\Sigma$ over $\Sigma$, $p=(p(0)p(1)\ldots)$.
For a word $w\in\s$ let $|w|$ be its length and let $\varepsilon\in\s$ be the empty word. For $p\in \om$ let $p^ {<i}\in\s$ be the prefix of $p$ of
length $i\in\IN$. We use the  ``wrapping function''
$\iota:\s\to\s$, $\iota(a_1a_2\ldots a_k):=110a_10a_20\ldots a_k011$
for coding words such that $\iota(u)$ and $\iota(v)$ cannot overlap properly.
Let $\langle i,j\rangle:=(i+j)(i+j+1)/2+j$ be the bijective Cantor pairing function on $\IN$.
We consider standard functions for finite or countable tupling on $\s$ and $\om$
denoted by $\langle\,\cdot\,\rangle\:$ \cite[Definition~2.1.7]{Wei00}, in particular,
$\langle u_1,\ldots, u_n\rangle=\iota(u_1)\ldots\iota(u_n)$, $\langle u,p\rangle=\iota(u)p$,
$\langle p,q\rangle= (p(0)q(0)p(1)q(1)\ldots)$ and
$\langle p_0,p_1,\ldots\rangle \langle i,j\rangle=p_i(j)$
 for $u,u_1,u_2,\ldots\in\s$ and  $p,q,p_0,p_1,\ldots\in\om$.
Consider $u\in\s$ and $w\in \s\cup\om$. Let $u\prf w$ iff $\iota(u)$ is a prefix of $w$,
$u\ll w$ iff $\iota(u)$ is a subword of $w$
and let $\widehat w$ be the longest subword $v\in 11\s11$ of $w$
(and the empty word if no such subword exists). Then for $u,w_1,w_2\in\s$,
$(u\ll w_1\vee u\ll w_2)\iff u\ll\widehat w_1\widehat w_2$.

For $Y_0,\ldots,Y_n\in\{\s,\om\}$ a partial function $f\pf Y_1\times\ldots\times Y_n\to Y_0$ is computable, if it is computed by a Type-2 machine. A Type-2 machine $M$ is a Turing machine with $n$ input tapes, one output tape and finitely many additional work tapes.
A specification assigns to the input tapes $1,\ldots ,n$ and the output tape $0$
types $Y_i\in\{\s,\om\}$ such that the machine computes a  function $f_M\pf Y_1\times\ldots\times Y_n\to Y_0$ \cite{Wei00}. Notice that on the output tape the machine can only write and move its head to the right.

A notation of a set $X$ is a surjective partial function $\nu\pf\s\to X$ and a representation is a surjective partial function $\delta\pf\om\to X$. Here, finite or infinite sequences of symbols are considered as ``concrete names'' of the ``abstract'' elements of $X$. Computability  on $X$ is defined by computations on names.
Let $\gamma_i\pf Y_i\to X_i$, $Y_i\in\{\s,\om\}$ for $i\in\{0,1\} $ be notations or representations. A set $W\In X_0$ is called $\gamma_0$-r.e. (recursively enumerable), if
there is a Type-2 machine $M$ that halts on input $y_0\in\dom(\gamma_0)$ iff $\gamma_0(y_0)\in W$.
A function $h\pf Y_1\to Y_0$ realizes  a multi-function $f:X_1\mto X_0$, iff $\gamma_0\circ h(y_1)\in f\circ \gamma_1(y_1)$ whenever
$f\circ \gamma_1(y_1))\neq \emptyset$. The function $f$ is called $(\gamma_1,\gamma_0)$-computable, if it has a computable realization. The definitions can be generalized straightforwardly to subsets of $X_1\times\ldots\times X_n$
 and multi-functions $f:X_1\times\ldots\times X_n\to X_0$ ($(\gamma_1,\ldots,\gamma_n)$-r.e., $(\gamma_1,\ldots,\gamma_n,\gamma_0)$-computable).

In this article we study axioms of computable separation for  {\em computable topological spaces}  ${\bf X}=(X,\tau,\beta,\nu)$  \cite[Definition~4]{WG09}, where $\tau$ is a $T_0$-topology on the set $X$ and $\nu\pf\s\to\beta$ is a notation of a base $\beta$ of $\tau$ such that $\dom(\nu)$ is recursive and there is an r.e. set $S\In (\dom(\nu))^3$ such that $\nu(u)\cap\nu(v)=\bigcup\{\nu(w)\mid (u,v,w)\in S\}$.
We mention expressly that in the past various spaces have been called ``computable topological space''.
We allow $U=\emptyset$ for $U\in\beta$ which is forbidden, for example, in \cite{GWX07a,XG09}.

We define a notation $\nufs$ of the finite subsets of the base $\beta$ by
$\nufs(w)=W$ $:\iff ((\forall v\ll w)v\in\dom(\nu) \an W=\{\nu(v)\mid v\ll w\})$.
Then $\bigcup\nufs$ and $\bigcap\nufs$ are notations of the finite unions and the finite intersections of base elements, respectively.

For the points of $X$ we consider the canonical (or inner) representation $\delta\pf\om\to X$; $\delta(p)=x$ iff $p$ is a list of all $\iota(u)$ (possibly padded with 1s) such that $x\in\nu(u)$ (hence $u\ll p\iff \delta(p)\in\nu(u)$) . For the set of open sets, the topology $\tau$ we consider the inner representation $\theta\pf\om\to\tau$  defined by $u\in\dom(\nu)$ if $u\ll p\in\dom(\theta)$ and $\theta(p):=\bigcup\{\nu(u)\mid u\ll p\}$. For the closed sets we consider the outer representation $\psi^-(p):=X\setminus \theta(p)$
\cite{WG09}.

The canonical notations of the natural and the rational numbers are denoted by $\nu_\IN$ and $\nu_\IQ$, respectively. For the real numbers we use the canonical representation $\rho$ (Example~\ref{e8}(\ref{e8a})), the lower representation $\rho_<$ and the upper representation $\rho_>$ \cite{Wei00}.

\section{\boldmath Axioms of Computable separation}\label{secf}

For a topological space ${\bf X}=(X,\tau)$ with set  ${\mathcal A}$ of closed sets we consider the following separation properties:

\begin{defi}[axioms of separation]\label{d50}
\begin{eqnarray*}
&{\rm T_0:}&(\forall x,y\in X,\ x\neq y) (\exists W\in\tau)
((x\in W\wedge y\not\in W)\vee (x\not\in W\wedge y\in W))),\\
&{\rm T_1:}& (\forall x,y\in X,\ x\neq y) (\exists W\in\tau)
(x\in W\wedge y\not\in W),\\
&{\rm T_2:}& (\forall x,y\in X,\ x\neq y) (\exists U,V\in\tau)
(U\cap V=\emptyset \wedge x\in U\wedge y\in V),\\
&{\rm T_3:}
&(\forall x\in X,\forall A\in{\mathcal A}, x\not\in A) (\exists U,V\in\tau)
(U\cap V=\emptyset \wedge x\in U\wedge A\In V),\\
&{\rm Ty:}
&(\forall x\in X,\forall A\in{\mathcal A}, x\not\in A) (\exists f:X\to\IR)\\
&&(\mbox{ $f$ is continuous,
$\range(f)\In [0;1]$, $f(x)=0$ and $f[A]\In\{1\}$})\\
&{\rm T_4: }
& (\forall  A,B\in{\mathcal A}, A\cap B=\emptyset) (\exists U,V\in\tau)
(U\cap V=\emptyset \wedge A\In U\wedge B\In V).\\
&{\rm Ur:}
&(\forall  A,B\in{\mathcal A}, A\cap B=\emptyset) (\exists f:X\to\IR)\\
&&(\mbox{ $f$ is continuous,
$\range(f)\In [0;1]$, $f[A]\In\{0\}$ and $f[B]\In\{1\}$.})
\end {eqnarray*}
We will speak of $\,T_2$-spaces, $\Ty$-spaces etc.
\end{defi}
$T_2$-spaces are called {\em Hausdorff spaces}.
 Many authors, for example \cite{Eng89}, call a space $T_3$-space or regular iff $\rm T_1+T_3$, call a space  $T_{3\frac{1}{2}}$-space, Tychonoff space or completely regular
iff $\rm T_1+Ty$, and call a space $T_4$-space or normal iff $\rm T_1+T_4$.
From topology \cite{Eng89} we know:
\[\rm T_1+Ur\iff T_1+T_4\Longrightarrow T_1+Ty\Longrightarrow T_1+T_3\Longrightarrow
\rm T_2\Longrightarrow T_1\Longrightarrow T_0\]
where the implications are proper. The first implication from the right to the left is
{\em Urysohn's lemma}.
We mention that $(X,\tau)$ is a $T_1$-space, iff all sets
$\{x\}$ ($x\in X$) are closed \cite{Eng89}.

In this article we consider only computable topological spaces ${\bf
X}=(X,\tau,\beta,\nu)$, which are $T_0$-spaces with countable base (also called {\em second countable}). For such spaces
$\rm T_3\Longrightarrow T_2$ and $\rm T_1 +T_4\iff T_1+T_3$ \cite[Theorem~1.5.16]{Eng89}, hence
\begin{eqnarray}\label{f86}
\rm T_1+Ur\iff T_1+T_4\iff T_1+Ty\iff T_3.
\end{eqnarray}

Axioms of computable separation for $T_0$, $T_1$ and $T_2$ have been studied in \cite{Wei10}. In the following we introduce computable versions of the axioms
$\rm T_3,\ T\hspace{-.3ex}y,\ T_4$ and $\rm Ur$. The  computable Hausdorff axioms
$\rm CT_2$ and $\rm SCT_2$ are from \cite{Wei10}.
In the direct effectivizations the existing objects must be computed. For the points we compute basic neighborhoods (w.l.o.g.) instead of general open neighborhoods.
Let $C(X,\IR)$ be the set of continuous functions $f:X\to\IR$ and let $[\delta\to\rho]$ be the canonical representation of this set \cite{WG09,Wei00}.

\begin{defi}[axioms of computable separation]\label{d4}$ $\\[1ex]
\noindent $\rm CT_2:$
The multi-function $t_2:X\times X\mto \beta\times \beta$ is $(\delta,\delta, [\nu,\nu])$-computable where
\\ \phantom{$\rm VT_2:$} $(U,V)\in t_2(x,y)$ iff $x\in U$, $y\in V$ and $U\cap V=\emptyset$.
\\[-1.5ex]

\noindent $\rm SCT_2:$   There is an r.e. set $H\In \s \times \s$ such that
\vspace{-1.5ex}
\begin{eqnarray}
\label{f69} &&
(\forall x,y,\ x\neq y)(\exists  (u,v)\in H)( x\in\nu(u)\wedge y\in\nu(v))\ \ \mbox{ and }\\
\label{f70}
&&(\forall  (u,v)\in H)\: \nu(u)\cap \nu(v)=\emptyset \,.
\end{eqnarray}

\vspace{-1.0ex}
\noindent $\rm WCT_3:$
The multi-function $t^w_3:X\times \beta\mto \beta$ is $(\delta,\nu, \nu)$-computable where
\\ \phantom{$\rm WCT_3:$} $U\in t^w_3(x,W)$ iff $x\in U\In \overline U\In W$.
\\[-1.5ex]

\noindent $\rm CT_3:$
The multi-function $t_3:X\times{\mathcal A}\mto \beta\times \tau$ is $(\delta,\psi^-, [\nu,\theta])$-computable,
\\ \phantom{$\rm CT_3:$} where $(U,V)\in t_3(x,A)$ iff  $\,x\not\in A$,
$U\cap V=\emptyset$, $x\in U$ and $A\In V$.
\\[-1.5ex]

\noindent $\rm CT_3':$
The multi-function $t'_3:X\times \beta\mto \beta\times {\mathcal A}$ is $(\delta,\nu, [\nu,\psi^-])$-computable \\
\phantom{$\rm SCT_3:$}where $(U,B)\in t'_3(x,W)$ iff  $\,x\in U\In B\In W$.
\\[-1.5ex]

\noindent $\rm SCT_3:$ There are an r.e. set $R\In \dom(\nu)\times \dom(\nu)$ and a computable function\\
\phantom{$\rm SCT_3:$} $r\pf \s\times\s\to\om$ such that for all $u,w\in\dom(\nu)$,
\vspace{-1ex}
\begin{eqnarray}
\label{f39}
&& \nu(w)=\bigcup\{\nu(u)\mid(u,w)\in R\}\,,\\
\label{f40}&&
 (u,w)\in R \Longrightarrow \ \nu(u)\In\psi^-\circ r(u,w)\In \nu(w)\,.
\end{eqnarray}

%\vspace{-1ex}
\noindent $\rm CTy:$ The multi-function $t_{\rm Ty}:X\times {\mathcal A}\mto C(X,\IR)$
is $(\delta,\psi^-,[\delta\to\rho])$-computable
\\ \phantom{$\rm CTy$: }where $f\in t_{\rm Ty}(x,A)$ iff $\range(f)\In[0;1]$, $\,x\not\in A$, $f(x)=0$ and $f[A]\In\{1\}$.\\[-1.5ex]

\noindent $\rm CTy'\hspace{-.4ex}:$ The multi-function
$t_{\rm Ty}':X\times \beta \mto \beta\times C(X,\IR)$
is $(\delta,\nu,[\nu,[\delta\to\rho]])$-computable
\\ \phantom{$\rm CTy'$: }where $(U,f)\in t_{\rm Ty}'(x,W)$ iff $\range(f)\In[0;1]$,
$\,x\in U\In W$, $f[U]=\{0\}$
\\ \phantom{$\rm CTy'$: }and $f[X\setminus W]\In\{1\}$.\\[-1.5ex]

\noindent $\rm SCTy:$
There are an r.e. set $T\In \dom(\nu)\times \dom(\nu)$ and a computable
\\\phantom{$\rm SCTy$: }function $t\pf \s\times\s\to\om$ such that
\vspace{-1ex}
\begin{eqnarray}
\label{f87}& \nu(w)=\bigcup\{\nu(u)\mid(u,w)\in T\}& \mbox{for all }\  w\in\dom(\nu) \ \ \mbox{and}\\
\label{f88} &
f_{uw}[\nu(u)]\In\{0\} \  \mbox{and}\  f_{uw}[X\setminus\nu(w)]\In\{1\} &\ \mbox{ for all }\  (u,w)\in T\,,
\end{eqnarray}
\phantom{$\rm SCTy$: } where $f_{uw}:=[\delta\to\rho]\circ t(u,w)$.
\\[-1.5ex]

\noindent $\rm CT_4:$
The multi-function $t_4:{\mathcal A\times A}\mto \tau\times\tau$ is $(\psi^-,\psi^-, [\theta,\theta])$-computable\\
\phantom{$\rm CT_4:$ }where $(U,V)\in t_4(A,B)$ iff
$\,U\cap V=\emptyset$, $A\In U$ and $ B\In V$.
\\[-1.5ex]

\noindent $\rm CUr:$ The multi-function $t_{\rm Ur}:{\mathcal A\times A}\mto C(X,\IR)$
is $\,(\psi^-,\psi^-,[\delta\to\rho])$-computable, \\
\phantom{$\rm CUr$: }where $f\in t_{\rm Ur}(A,B)$ iff $\range(f)\In[0;1]$,
$ A\cap B=\emptyset$, $f[A]\In\{0\}$
and $f[B]\In\{1\}$.
\end{defi}

The axioms $\rm CT_2$, $\rm CT_3$, $\rm CTy$, $\rm CT_4$ and $\rm CUr$ are the direct effectivizations of $T_2$, $\rm T_3$, $\rm Ty$, $\rm T_4$ and $\rm Ur$, respectively. Obviously, $\rm SCT_2$ implies $\rm T_2$. \  $\rm WCT_3$, $\rm CT_3 $, $\rm CT_3'$ and  $\rm SCT_3$ imply   $\rm T_3$. \
$\rm CTy$, $\rm CTy' $ and $\rm SCTy$ imply   $\rm Ty$. \  $\rm CT_4$ implies $\rm T_4$. \ $\rm CUr$ implies $\rm Ur$.
In contrast to ${\rm CT}'_3$, in ${\rm WCT}_3$ the function $t_3^w$ does not compute a
$\psi^-$-name of a closed set such that $x\in U\In B\In W$.
The sets $H$ from $\rm SCT_2$, $R$ from $\rm SCT_3$ and $T$ from $\rm SCTy$  may contain pairs $(u,w)$ such that $\nu(u)=\emptyset$ or $\nu(w)=\emptyset$.
Also, empty open or closed sets are not excluded  as inputs for the separating functions.

We do not consider the numerous variants of the separation axioms  where in some places the representations $\delta$ of the points,
$\theta$ of the open sets and $\psi^-$ of the closed sets are replaced by
$\delta^-$, $\theta^-$ and $\psi^+$, respectively \cite[Definition~5]{WG09}.
The following examples illustrate the definitions. Further examples are given in Section~\ref{secp}.

\begin{exa}\label{e8}\rm \hfill
\begin{enumerate}[(1)]
\item\label{e8a} The {\em computable real line} is defined by
${\bf R}:=(\IR,\tau_\IR,\beta,\nu)$ such that $\tau_\IR$ is the real
line topology and $\nu$ is a canonical notation of the set of all
open intervals with rational endpoints. ${\bf R}$  is a computable topological space.
Its canonical representation is called~$\rho$. All the axioms from Definition~\ref{d4} are true for~${\bf R}$.

\item\label{e8e}  A computable metric space is a tuple ${\bf M}=(X,d,A,\alpha)$ such that $(M,d)$ is a metric space and $\alpha$ is a notation with recursive domain of a set $A$ which is dense in $X$ such that the distance d restricted to $M\times M$ is
$(\alpha,\alpha,\rho)$-computable  \cite[Definition~8.1.2]{Wei00}. Let $\nu$ be a
canonical notation of the set $\beta$ of all open balls with center from $A$ and rational radius and let $\tau$ be the smallest topology containing $\beta$.
Then ${\bf X}=(X,\tau,\beta,\nu)$ is a computable topological space for which all the axioms from Definition~\ref{d4} are true (Theorem~\ref{t3}).

\item\label{e8b} ($CT_0$ and $CT_4$ but not $T_1$, $T_2$ or $T_3$) A space is $ CT_0$
iff the multi-function $t_0$ is $(\delta,\delta,\nu)$-computable, where
$t_0$ maps every $(x,y)\in X^2$ such that $ x\neq y$ to some
$U\in \beta$ such that ($x\in U$ and $y\not\in U)$ or ($x\not\in U$ and  $y\in U$)
\cite{Wei10}.

Let ${\bf Si}:=(\{\bot,\top\},\tau_{\bf Si},\beta_{\bf Si},\nu_{\bf Si})$ be
the Sierpinski space defined by $\nu_{\bf Si}(0)=\{\bot,\top\}$ and $\nu_{\bf Si}(1)=\{\top\}$. The space  is $T_0$ but not $T_1$.

There is a machine $M$ that on input $(p,q)\in\om\times\om$ writes $1$ and halts.
The function $f_M$ realizes the function $(x,y)\mapsto \{\top\}$. Then for $x\neq y$, ($x=\top$ and $y=\bot$) or ($x=\bot$ and $y=\top$), hence for $U:=\{\top\}=\nu(1)$,
($x\in U$ and $y\not\in U)$ or ($x\not\in U$ and  $y\in U$). Therefore, $\bf Si$ is $CT_0$.

There are computable sequences $p',q'\in\om$ such that $\theta_{\bf Si}(p')=\emptyset$ and $\theta_{\bf Si}(q')=\{\bot,\top\}$. There is a machine $M$
that on input $(p,q)$ searches in $p$ and $q$ until it has found $0\ll p$ or $0\ll q$.
In the first case it writes $\langle p',q'\rangle$ and in the second case $\langle q',p'\rangle$.
Let $\psi^-_{\bf Si}(p)=A$ and $\psi^-_{\bf Si}(q)=B$ such that $A\cap B=\emptyset$. Then $A=\emptyset$ or $B=\emptyset$, hence $0\ll p$ or $0 \ll q$. In the first case,
$A=\emptyset\In \theta_{\bf Si}(p')$ and $B\In \theta_{\bf Si}(q')=\{\bot,\top\}$ and in the second case, $A\In \theta_{\bf Si}(q')=\{\bot,\top\}$ and $B=\emptyset\In \theta_{\bf Si}(p')$.
Therefore,  $f_M$ realizes~$t_4$.

\item\label{e8c} (discrete implies $WCT_3$) Let $\bf X$ be a discrete computable topological space. Then every subset of $X$ is open and closed and $X$ is countable. The function $t^w_3:(x,W)\mapsto W$ is $(\delta,\nu,\nu)$-computable. It satisfies $\rm WCT_3$ since $x\in W$ implies
$x\in\overline W\In W$.

\end{enumerate}
\end{exa}

\noindent By the next lemma the above computable separation axioms are robust, that is, they do not depend on the notation
$\nu$ of the base explicitly but only on the computability concept on the points induced by it. Two computable topological spaces  ${\bf X}=(X,\tau,\beta,\nu)$ and
$\widetilde{\bf X}=(\widetilde X,\widetilde\tau,\widetilde\beta,\widetilde\nu)$ are called equivalent  iff $(X,\tau)=(\widetilde X,\widetilde\tau)$, $\nu\leq \widetilde \theta$ and $\widetilde\nu\leq \theta$, that is, there are  computable functions
$g,\widetilde g\pf\s\to\om$ such that
\begin{eqnarray}\label{f23}
\nu(u)=\widetilde\theta\circ g(u)&\mbox{and}& \widetilde\nu(u)=\theta\circ \widetilde g(u)\,.
\end{eqnarray}
The condition ``$\nu\leq \widetilde \theta$ and $\widetilde\nu\leq \theta$''
is equivalent to $\delta\equiv \widetilde\delta$.
For equivalent topological spaces,
$\theta\equiv \widetilde\theta$, $\psi^-\equiv \widetilde\psi^-$ and
 $\kappa\equiv \widetilde\kappa$ \cite[Definition~21, Theorem~22]{WG09}.

\begin{lem}\label{l5}
 Let $\widetilde{\bf X}=(X,\tau,\widetilde\beta,\widetilde\nu)$ be a computable topological space
equivalent  to ${\bf X}=(X,\tau,\beta,\nu)$ .
Then each separation axiom from Definition~\ref{d4} for ${\bf X}$ is equivalent to the corresponding axiom for $\widetilde{\bf X}$.
\end{lem}

\pproof

$\bf SCT_2$: See \cite{Wei10}.

$\bf WCT_3$: Assume $\rm \widetilde{WCT}_3$. Let $x=\delta(p)$, $W=\nu(w)$ and $x\in W$. Since $\delta\equiv \widetilde\delta$ and $\nu\leq \widetilde \theta$ we can compute some $\widetilde p$ and some $\widetilde w$ such that
$x=\widetilde \delta(\widetilde p)\in\widetilde\nu(\widetilde w)\In\nu(w)$.
By $\rm \widetilde{WCT}_3$ we can compute some $\widetilde u$ such that
$x\in \widetilde \nu(\widetilde u)\In {\rm closure}(\widetilde \nu(\widetilde u)) \In \widetilde\nu(\widetilde w)$. Since $\widetilde \nu\leq \theta$, from $p$ and $\widetilde u$ we can compute some $u$ such that $x\in\nu(u)\In \nu(\widetilde u)$.
We obtain $x\in\nu(u)\In\overline{\nu(u)}\In\nu(w)$. Therefore, $\rm WCT_3$ is true.
By symmetry, $\rm WCT_3\Longrightarrow \widetilde{WCT}_3$\,.

$\bf SCT_3$:
Assume $\rm SCT_3$. With the functions $g, \widetilde g$ from (\ref{f23}) let
$\widetilde R:=\{(\widetilde u,\widetilde w)\mid (\exists (u,w)\in R)(w\ll \widetilde g(\widetilde w),\  \widetilde u\ll g(u))\}$. Then $\widetilde R$ is r.e.
Suppose $(\widetilde u,\widetilde w)\in\widetilde R$. Then for some $(u,w)\in R$,
$\widetilde\nu(\widetilde u)\In \widetilde \theta \circ g(u) =\nu(u)\In\nu(w) \In
\theta\circ \widetilde g(\widetilde w)=\widetilde\nu(\widetilde w)  $.
On the other hand suppose, $x\in\widetilde\nu(\widetilde w)=\theta\circ\widetilde g(\widetilde w)$. Then $x\in\nu(w)$ for some $w\ll \widetilde g(\widetilde w)$.
By ${\rm SCT}_3$ there is some $u$ such that $(u,w)\in R$ and $x\in\nu(u)=\widetilde \theta\circ g(u)$. Then $x\in\widetilde \nu(\widetilde u)$ for some $\widetilde u\ll g(u)$. In summary, $x\in\widetilde \nu(\widetilde u)$ for some $\widetilde u$ such that
$(\widetilde u,\widetilde w)\in\widetilde R$.
Therefore,  (\ref{f39}) holds for $\widetilde \nu$ and $\widetilde R$.

There is a computable function $d$ translating $\psi^-$ to $\widetilde \psi^-$ \cite{WG09}.
Let $M$ be a machine that on input $(\widetilde u,\widetilde w)$ searches for $(u,w)\in R$ such that $w\ll \widetilde g(\widetilde w)$ and  $\widetilde u\ll g(u)$ and then computes $d\circ r(u,w)$. Then
$\widetilde\nu(\widetilde u)\In \widetilde \theta\circ g(u )=\nu(u)\In\psi^-\circ r(u,w)
\In\nu(w)\In \theta\circ \widetilde g(\widetilde w)=\widetilde \nu(\widetilde w)$. Since $\psi^-\circ r(u,w)= \widetilde \psi^-\circ  d \circ r(u,w)= \widetilde \psi^-\circ
f_M(\widetilde u,\widetilde w)$,
(\ref{f40}) holds for $\widetilde {SCT}_3$ with $\widetilde r:= f_M$ and $\widetilde R$.
Therefore, $\rm SCT_3\Longrightarrow \rm\widetilde {SCT}_3$.
By symmetry, $\rm\widetilde {SCT}_3\Longrightarrow \rm SCT_3$.

For the other axioms the proofs are similar.
Notice that $[\delta\to\rho]\equiv [\widetilde\delta\to \rho]$ if
$\delta\equiv \widetilde\delta$.
\qq

\section{Implications}\label{seco}
In this section we prove a number of implications between the separation properties, in  Section~\ref{secp} we prove by counterexamples that some of the implications are proper.
A topological space is discrete iff every singleton $\{x\}$ is open iff every subset $B\In X$ is open. A discrete space is $T_i$ for $i=1,\ldots, 4$. Let $\rm D$ be the axiom stating that the space is discrete.

\begin{thm}\label{t6} $ $
\begin{enumerate}[\em(1)]
\item\label{t6a}
$\rm SCT_3\Longrightarrow CTy\Longrightarrow  CT_3\Longrightarrow SCT_2 \Longrightarrow CT_2$,
\item\label{t6f}
$\rm D \Longrightarrow WCT_3$,
\item\label{t6b}
$\rm CT_3\Longrightarrow WCT_3$
\item\label{t6c}
$\rm SCT_3\Longrightarrow CT_4$,
\item\label{t6d}
$\rm SCT_3\iff SCTy$,  \  $\rm CTy\iff CTy'$, \  $\rm CT_3\iff CT_3'$,
\item \label{t6e}$\rm CT_4\iff CUr$.
\end{enumerate}
\end{thm}

The implications $\rm SCT_3\Longrightarrow CT_4\Longrightarrow CUr$ have been already been proved in \cite{Sch98} for a computable topological space $T({\bf Z})$
derived from a predicate space $\bf Z$ (in the terminology of \cite{WG09}). For our computable topological space ${\bf X}=(X,\tau,\beta,\nu)$, ${\bf Z}:=(X,\beta,\nu)$ is a predicate space and $T({\bf Z})=(X,\tau,\tilde\beta, \tilde\nu)$, where $\tilde\nu$ is the notation of the finite intersections of base elements canonically derived from $\nu$, is equivalent to $\bf X$ by \cite[Lemma~23]{WG09}. By Lemma~\ref{l5},
$\rm SCT_3\Longrightarrow CT_4\Longrightarrow CUr$ for a computable topological space follows from \cite{Sch98}.
More concise proofs are given in \cite{GSW07} for a computable topological space
${\bf X}=(X,\tau,\beta,\nu)$ such that $U\neq \emptyset$ for all $U\in\beta$.
This restriction, however, is unnecessary. The reader may check this in  Appendix A.\\

\pproof

$\bf SCT_3 \Longrightarrow CT_4$: (cf. \cite[Lemma~1.5.15, Theorem~1.5.17]{Eng89})
The proof from \cite{GSW07} is added in Appendix B.

\medskip
$\bf CT_4\Longrightarrow CUr:$ See \cite[Theorem~1.5.15]{Eng89}.
The proof from \cite{GSW07} is added in the appendix.

\medskip
$\bf CUr\Longrightarrow CT_4$:
By the multi-function $t_{\rm UR}$ from $A,B$ such that $A\cap B=\emptyset$ we can compute a continuous function $f:X\to\IR$ such that
$\range(f)\In[0;1]$, $f[A]\In\{0\}$  and $f[B]\In\{1\}$.
Then by \cite[Theorem~38]{WG09} the open sets $U:=f^{-1}[(-\infty;1/2)]$ and $V:=f^{-1}[(1/2);\infty)]$ can be computed.
They separate $A$ and $B$.

\medskip
$\bf CUr + SCT_3\Longrightarrow SCTy$: Let $R$ be the set and let $r$ be the function from ${\rm SCT}_3$. Define $T:=R$. By $r$ from $(u,w)\in T$ we can compute a closed set
$A$ such that $\nu(u)\In A\In \nu(w)$. Since  $U\mapsto U^c$ for base sets is $(\nu,\psi^-)$-computable, by $t_{Ur}$ from $A$ and $X\setminus \nu(w)$ we can compute some continuous function $f:X\to \IR$ such that $\range(f)\In[0;1]$,
 $f[A]\In\{0\}$, hence $f[\nu(u)]\In \{0\}$  and $f[X\setminus \nu(w)]\In\{1\}$.

\medskip
$\bf SCTy\Longrightarrow SCT_3$:
Let $T$ be the set and $t$ be the function from $\rm SCTy$.
Define $R:= T$. By \cite[Theorem~38]{WG09} the function $f\mapsto f^{-1}(1/2,\infty)$ for continuous $f:X\to\IR$ is $([\delta\to\rho],\theta)$-computable. Let $h\pf \om\to\om$ be a computable realization.
Then  for $(u,w)\in R$ and $V:= ([\delta\to\rho]\circ t(u,w))^{-1}(1/2;\infty)
=\theta\circ h\circ t(u,w)$,
$\nu(u)\cap V=\emptyset $ and $X\setminus \nu(w)\In V$.
Therefore, $\nu(u)\In X\setminus V=\psi^-\circ h\circ t(u,w)=X\setminus V\In \nu(w)$. Define $r:= h\circ t$.

$\bf CTy' \Longrightarrow CTy$: From $(x,A)$ such that $x\not\in A$ some $W\in\beta$ can be
computed such that $x\in W\In X\setminus A$.
From $(x,W)$ such that $x\in W$  by $t_{\rm Ty}'$ some $U\in\beta$ and some
continuous function $f:X\to\IR$ can be computed such that $x\in U\In W$, $\range(f)\In [0;1]$,
$f(y)=0$ for $y\in U$ and $f(y)=1$ for $y\not\in W$. For this function $f$,
$f(x)=0$ since $x\in U$ and $f(y)=1$ for $y\in A$ since $W\In X\setminus A$.
Therefore, $t_{\rm Ty}$ is $(\delta,\psi^-,[\delta\to\rho])$-computable.

\medskip
$\bf CTy \Longrightarrow CTy'$: Suppose $x\in W\in \beta$. From $W$, $A:= X\setminus W$ can be computed.
From $(x,A)$ by $t_{\rm Ty}$ some continuous function $g$ can be computed such that $\range(g)\In [0;1]$, $g(x)=0$ and $g(y)=1$ for $y\in A$. Let $f(y):=\max(0,2g(y)-1)$. Then $g\mapsto f$ is
$([\delta\to\rho],[\delta\to\rho])$-computable. Obviously $\range(f)\In [0;1]$, $f(y)=0$ for $g(y)<1/2$ and $f(y)=1$ for $y\in A$. By \cite[Theorem~38]{WG09}, $g\mapsto g^{-1}(-\infty,1/2)$ is
$([\delta\to\rho],\theta)$-computable. Finally from $(x,V)$ such that $x\in V\in\tau$ some $U\in\beta$
can be computed such that $x\in U\In V$. Notice that $f(y)=0$ for $y\in U\In g^{-1}(-\infty,1/2)$. Therefore, from $(x,W)$ some $(U,f)\in t_{\rm Ty}(x,W)$ can be computed.

\medskip
{\boldmath $\rm CT_3' \Longrightarrow  CT_3$:} Using $t'_3$ from $(x,A)$ we can compute
in turn $(x,X\setminus A)$,  $(x,W)$ for some $W\in\beta$ such that $x\in W\In X\setminus A$, some
$(U,B)\in\beta\times{\mathcal A}$ such that $x\in U\In B\In W$, and finally $(U,V)$ where $V:=X\setminus B$.
By simple transformations, $U\cap V=\emptyset$, $x\in U$ and $A\In  V$.

\medskip
{\boldmath $\rm CT_3 \Longrightarrow  CT_3'$:} The function $W\mapsto X\setminus W$
is $(\nu,\psi^-)$-computable, and the function $V\mapsto X\setminus V$ is $(\theta,\psi^-)$-computable.
Then using $t_3$, from $(x,W)$ we can compute in turn $(x,A)$, $A:=X\setminus W$,
$(U,V)\in\beta\times \tau$ such that $U\cap V=\emptyset$, $x\in U$ and $A\In  V$, and $(U,B)$, $B:=X\setminus V$. By simple transformations, $x\in U\In B\In W$.

\medskip
{\boldmath $\rm CT_3' \Longrightarrow  WCT_3$:} Obvious.

\medskip
$\bf SCTy \Longrightarrow CTy'$: Let $T$ be the set and $t$ be the function from ${\rm SCTy}$. Assume
$x=\delta(p)\in \nu(w)$. By (\ref{f87}) there is some $u\in\dom(\nu)$ such that $u\ll p$,
and $(u,w)\in T$.
Then $f_{uw}:=[\delta\to\rho]\circ t(u,w)$ satisfies (\ref{f88}).
There is a machine that on input $(p,w)$ searches for some $u$ such that $u\ll p$ and
$(u,w)\in T$ and writes $\langle u,t(u,w)\rangle$. Then $f_M$ realizes $t_{\rm Ty}'$.

\medskip
$\bf CTy \Longrightarrow CT_3$: Suppose $x\not\in A$ and $A$ is closed. By $t_{\rm Ty}$ we can compute
some continuous function $f$ such that $f(x)=0$ and $f(y)=1$ for $x\in A$. By \cite[Theorem~38]{WG09}, the functions  $f\mapsto f^{-1}(-\infty,1/2)$ and $f\mapsto f^{-1}(1/2,\infty)$ are $([\delta\to\rho],\theta)$-computable. Since $x\In f^{-1}(-\infty,1/2)$ and $A\In f^{-1}(1/2,\infty)$,
the multi-function $(x,A)\mmto (U,V)$ such that $x\in U$, $A\In V$  and $U\cap V=\emptyset$ is
$(\delta,\psi^-,[\theta,\theta])$-computable. From $q$ and $r$ such that $\delta(q)\in\theta(r)$ we can compute some $u$ such that $x\in\nu(u)\In \theta(r)$.
Therefore, $t_3$ is
$(\delta,\psi^-,[\nu,\theta])$-computable.

\medskip
{\boldmath $\rm CT'_3 \Longrightarrow  SCT_2$:} Let $M$ be a machine such that $f_M$ realizes~$t'_3$.
Since finite intersection is $(\nufs,\theta)$-com\-putable \cite[Theorem~11]{WG09},
there is a computable function $g$ such that $\nui(w)=\theta\circ g(w)$.
Let $H$ be the set of all $(u,v)\in\dom(\nu)\times \dom(\nu)$ with the following properties: there are words $w,u_1,v_1,v_2$ such that
$v_1\in\dom(\nufs)$, $w\ll v_1$, $u \ll g(v_1)$ and on input $(v_11^\omega,w)$ in ${\rm length}(v_1)$ steps the machine $M$ writes at least $\iota(u_1)v_2$ such that $u_1 \ll v_1$, and $v\ll v_2$. The set $H$ is r.e.

Suppose $\delta(p)=x\neq y$.
Since $\rm CT_3'\Longrightarrow T_3$ and $\rm T_3\Longrightarrow T_2$ for second countable spaces, the space is $T_2$, hence there  is some $w$ such that $x\in\nu(w)$
and $y\not \in\nu(w)$. Then on input $(p,w)$ the machine $M$ writes some $\iota(u_1)q$
such that $x\in\nu(u_1)\In \psi^-(q)\In \nu(w)$.
Since $y\not\in\nu(w)$, hence $y\in\theta(q)$, there are a prefix $v_2$ of $q$ and a word $v$ such that $v\ll v_2$ and $y\in\nu(v)$. For producing
$\iota(u_1)v_2$ some prefix of $p$ is sufficient.
Since $x\in\nu(u_1)\In \nu(w)$ there is a prefix $v_1$ of $p$ such that $w\ll v_1$, $u_1 \ll v_1$ and  on input $(v_11^\omega,w)$
in ${\rm length}(v_1)$ steps the machine $M$ writes at least $\iota(u_1)v_2$.
Since $x\in\theta\circ g(v_1)$, there is some $u\ll g(v_1)$ such that $x\in \nu(u)$.
By definition of $H$, $(u,v)\in H$, hence (\ref{f69}) is true.

Suppose $(u,v)\in H$. Then there are words $w,u_1,v_1,v_2$ with the properties listed in the definition of~$H$. If $\nu(u)=\emptyset$, (\ref{f70}) is true.

Suppose $x\in\nu(u)\neq\emptyset$.
Since $u \ll g(v_1)$ and $w\ll v_1$,
$x\in\nu(u)\In  \theta\circ g(v_1) =\nui(v_1)\In \nu(w)$.
There is some $p'\in\om$ such that
$x=\delta(v_1p')\in \nui(v_1)$.
Since $M$ realizes $t'_3$,
on input $(v_1p',w)$  the machine $M$ writes $\iota(u_1)q$ such that
$x\in\nu(u_1)\In\psi^-(q)\In\nu(w)$.
In ${\rm length}(v_1)$ steps the machine can read only symbols from $v_1$ and, therefore, has the same behavior on input $(v_11^\omega,w)$.
By assumption on $u,v,w,u_1,v_1$ and $v_2$,
$v_2$ is a prefix of $q$ and $u_1\ll v_1$,
hence $\nu(u)\In\nu(u_1)$.
Since $v\ll v_2$, $v\ll q$, hence $\nu(v)\In \theta(q)$.
Since $\nu(u)\In\nu(u_1)\In \psi^-(q)$,
$\nu(u)\cap\nu(v)\In \nu(u_1)\cap\theta(q)=\emptyset$. Therefore,  (\ref{f70}) is true.

\medskip
$\bf SCT_2\Longrightarrow CT_2$: See \cite{Wei10}

\medskip
$\bf D\Longrightarrow WCT_3$: See Example~\ref{e8}(\ref{e8c}).

\medskip
\noindent The statements (\ref{t6a}) -- (\ref{t6e}) of the theorem follow from these results.
\qq\\

By \cite[Theorem~7]{Wei10}, the following statements are equivalent: ${\bf X}$ is $SCT_2$; $x\neq y$ is $(\delta,\delta)$-r.e.; $x\mapsto \overline{\{x\}}$ is $(\delta,\psi^-)$-computable. We apply this result in the next proof.

\begin{thm}\label{t9}
${\rm CT_4 + SCT_2 \Rightarrow CT_3}$.
\end{thm}

\pproof Since the space is $SCT_2$ it is $T_1$ \cite[Theorem~5]{Wei10}, hence $\{x\}=\overline{\{x\}}$ for every point $x$. Therefore by the above characterization, from $x$ and $A$ such that $x\not\in A$  we can compute $\{x\}$ and $A$ and by $CT_4$ we can compute disjoint open sets $U,V$ such that $\{x\}\In U$, hence $x\in U$, and $A\In V$. Therefore, the space is $CT_3$.
\qq\\

By \cite[Theorem~7]{Wei10}, $ {\rm T_2 \,\Rightarrow \, SCT_2}$ if $U\cap V=\emptyset$ is $(\nu,\nu)$-r.e.. A similar result holds for $T_3$-spaces.

\begin{thm}\label{t10}
 ${\rm CT_3\iff WCT_3}$ if $U\cap V=\emptyset$ is $(\nu,\nu)$-r.e.
\end{thm}

\pproof
Suppose ${\rm WCT_3}$. Then from $x\in X $ and $W\in \beta$ such that $x\in W$ we can compute some $U\in\beta$ such that $x\in U\in \overline U\in W$. For showing ${\rm CT_3'}$ it suffices to find a $\psi^-$-name of~$\overline U$. By assumption, from $U\in\beta$ we can find a list (encoded by $q\in\om$) of all $V\in\beta$ such that $U\cap V =\emptyset$.
Since for open $V$, $U\cap V=\emptyset \iff \overline U\cap V =\emptyset$, $q$ is a $\psi^-$-name of the closed $\overline U$. Therefore, the space is $CT_3'$, hence $CT_3$.
\qq\\

For a computable topological space ${\bf X}=(X,\tau,\beta,\nu)$ possibly $U=\emptyset$ for some $U\in\beta$.

\begin{thm}\label{t7}
 If the set $\{w\in\s\mid \nu(w)\neq\emptyset\}$ is r.e. then
 \[{\rm CT_3\iff CTy\iff SCT_3}\,.\]
\end{thm}

In particular, if all base elements are not empty then $\rm CT_3\iff SCTy\iff SCT_3$.
Of course, the space $\bf X$ in Example~\ref{e14} has empty base elements. The non-empty ones are not even r.e.\\

\pproof
Suppose that $\{ w\mid \nu(w)\neq\emptyset\}$ is r.e. Since finite intersection is $(\nufs,\theta)$-com\-putable \cite[Theorem~11]{WG09},
there is a computable function $g$ such that $\nui(w)=\theta\circ g(w)$.
Therefore, the set $\{w\in \s\mid \nui(w)\neq \emptyset\}$ is r.e.
Suppose the space is $CT_3$. By Theorem~\ref{t6} it is $CT'_3$.
There is a machine $M$ such that $f_M$ realizes the multi-function $t'_3$ from ${\rm CT}'_3$ in Definition~\ref{d4}.

Let $x_0=\delta(p_0)\in\nu(w)$. Then for some $u\in\dom(\nu)$ and $q_{p_0}\in\dom(\psi^-)$,
$f_M(p_0,w)=\langle u,q_{p_0}\rangle=\iota(u)q_{p_0}$ such that
\begin{eqnarray}
\label{f41}
x_0\in \nu(u)\In \psi^-(q_{p_0})\In \nu(w)\,.
\end{eqnarray}
For computing $\iota(u)$ some prefix $u_0\in\dom(\nufs)\cap \s 11$ of $p_0$ suffices. Since $\delta(p_0)\in \nu(w)$ we may assume $w\ll u_0$.
For all $p\in I_0:=\{p\in\dom(\delta)\mid u_0 \mbox{ is a prefix of } p\}$,
$f_M(p,w)=\iota(u)q_p$ for some $q_p$ such that
$\delta(p)\in \nu(u)\In \psi^-(q_p)\In \nu(w)$. Then
\begin{eqnarray}
x_0\in\nu(u)\in\bigcap_{p\in I_0}\psi^-(q_p)\In \nu(w)\,.
\end{eqnarray}
A word $u_0\in \dom(\nufs)$ is a prefix of some $p\in\dom(\delta)$ iff $\;\bigcap \nufs(u_0)\neq  \emptyset$.
We will define $R$ such that $(u,w)\in R$ iff for some $u_0\in \s 11$ such that $\;\bigcap \nufs(u_0)\neq  \emptyset$ the machine $M$ on input $(u_01^\omega,w)$ writes $\iota(u)$ in at most $|u_0|$ steps. From this word $u_0$ we will compute a sequence $q\in\om$ such that
$\psi^-(q)=\bigcap_{p\in I_0}\psi^-(q_p)$.

There is a machine $N$ that works on input $(u,w)$ as follows:\\
(S1) $N$ searches for some $u_0\in \dom(\nufs)\cap\s11$ such that $w\ll u_0$,  $\;\bigcap \nufs(u_0)\neq  \emptyset$  and the machine $M$ on input $(u_01^\omega,w)$ writes $\iota(u)$ in at most $|u_0|$ steps. \\
(S2) Then $N$ writes every $\iota(v)$ such that there are words $u'$ and $v'$ such that\\
(S2a)  $u'\in \dom(\nufs)\cap\s11$, $\;\bigcap \nufs(u')\neq  \emptyset$  and $u_0\pref u'$  and\\
(S2b) $M$ on input $(u'1^\omega,w)$ in $|u'|$ steps writes $\iota(u)v'$ such that $v\ll v'$. \\
Furthermore, $N$ writes $11$ repeatedly in order to produce an infinite sequence if only finitely many words $v$ can be found. If no $u_0$ can be found the machine does not halt and writes nothing.
Let $r:=f_N$ and $R:=\dom(f_N)$. Then $R\In\dom(\nu)\times \dom(\nu)$ and $R$ is r.e.
We must prove (\ref{f39}) and (\ref{f40}).

We show (\ref{f39}): Suppose $x=\delta(p)\in\nu(w)$. Then for some $u$ and $q$,$f_M(p,w)=\iota(u)q$ such that
$x\in\nu(u)\In\psi^-(q)\In \nu(w)$. There is a prefix $u_0\in \s11$ of $p$ such that $w\ll u_0$ and
$M$ on input $(u_01^\omega,w)$ writes $\iota(u)$ in at most $|u_0|$ steps, hence $(u,w)\in\dom(f_N)=R$. Therefore, $x\in\nu(u)$ for some $u$ with $(u,w)\in R$.
We conclude $\nu(w)\In\bigcup\{\nu(u)\mid (u,w)\in R\}$.

On the other hand, let $(u,w)\in R$. Then there is some
$u_0\in \dom(\nufs)\cap\s11$ such that $\;\bigcap \nufs(u_0)\neq  \emptyset$  and the machine $M$ on input $(u_01^\omega,w)$ writes $\iota(u)$ in at most $|u_0|$ steps.
There is some $p'$ such that $u_0p'\in\dom(\delta)$. Then $f_M(u_0p',w)=\iota(u)q'$ for some $q'$ hence $\delta(u_0p')\in\nu(u)\In\nu(w)$. Therefore,
$\bigcup\{\nu(u)\mid (u,w)\in R\}\In\nu(w)$.\\
Combining the two results we obtain  (\ref{f39}).

We show (\ref{f40}): Suppose $(u,w)\in R=\dom(f_N)$ is the input of the machine $N$ and let $q:=f_N(u,w)$. First, $N$ finds some $u_0$ with the properties listed in (S1).

 Suppose, later $N$ writes $\iota(v)$ as described in (S2). Then there are words $u',v'$
and  a sequence $p'\in \om$ such that $u'p'\in\dom(\delta)$ and $M$ on input $(u'p',w)$ in at most $|u'|$ steps writes $v'$ such that $v\ll v'\prf q$ and
$\delta(u'p')\in\nu(u)\In  \nu(w)$, hence $\nu(u)\cap\nu(v)=\emptyset$.
Therefore, $\nu(u)\cap\nu(v)=\emptyset$ for all $v$ such that $v\ll q$. We obtain
$\nu(u)\In \psi^-(q)= \psi^-\circ  f_N(u,w)$.

There are some $p',q'$ such that $u_0p'\in\dom(\delta) $ and $M$ on input $(u_0p',w)$
writes $\iota(u)q'$ such that $\delta(u_0p')\in\nu(u)\In \psi^-(q')\In\nu(w)$.
Suppose, $v\ll q'$. Then there are words $u',v'$ such that the conditions (S2a) and (S2b) are satisfied, hence $v\ll f_N(u,w)$. Therefore,
$\psi^-\circ f_N(u,w)\In \psi^-(q')\In\nu(w)$.\\
Combining the results we obtain $\nu(u)\In \psi^-\circ  f_N(u,w)\In\nu(w)$. Therefore we have proved (\ref{f40}).
\qq\\

Notice that the proof works correctly since in (S1) we have guaranteed $\bigcap\nufs(u_0)\neq \emptyset$ hence $u_0\prf p$ for some $p\in\dom(\delta)$. The realization $f_M$ of $t'_3$ may give unreasonable results on $(p,w)$ if $p\not\in\dom(\delta)$.

\section{Counterexamples and Summary}\label{secp}

We show by counterexamples that some of the implications
from Theorem~\ref{t6} are proper. In \cite{Wei10} a $CT_2$-space is given that is not $SCT_2$, hence  $\rm SCT_2 \Longrightarrow CT_2$ is proper.

\begin{exa}\label{e4}\rm($SCT_2$ but not $T_3$)
We extend  \cite[Example~1.5.7]{Eng89}.
Let ${\bf R}=(\IR,\tau_\IR,\beta_\IR,\nu_\IR)$ be the computable real line from Example~\ref{e8}(\ref{e8a}). Let $S:=\{1/i\mid i\in\IZ, \ i\neq 0\}$,
$\sigma:=\beta_\IR\cup \{(-1\,;\,1)\setminus S\}$ with canonical notation $\lambda$.
Then ${\bf Z}=(\IR,\sigma,\lambda)$ is a computable predicate space
and $T({\bf Z})=:(\IR,\tau,\beta,\nu)$ is a computable topological space \cite[Definition~8, Lemma~9]{WG09}. Since $\beta_\IR$ is a subset of the topology generated by $\sigma$ and $\nu_\IR\leq \lambda$, ``$x\neq y$'' is
$(\delta_{\bf Z},\delta_{\bf Z})$-r.e., hence $(\delta,\delta)$-r.e. By \cite[Theorem~7.2]{Wei10},
$T({\bf Z})$ is a $SCT_2$-space. The space is not $T_3$ since
the point $0$ cannot be separated from the closed set $S$ by disjoint open sets since
$U\cap S\neq\emptyset$ for every neighborhood $U$ of $0$.
\qq
\end{exa}

Since the above $SCT_2$-space is not $T_3$ it is not $WCT_3$, $CT_3$ or $SCT_3$.
The space from Example~\ref{e3} below is $T_4$ and $SCT_2$, but not $WCT_3$.
First we prove a lemma. Let us call a function $f\pf \IQ\to\IQ$ a {\em lower separation function} for a real number $x>0$ if  $f$ is computable
(precisely, $(\nu_\IQ,\nu_\IQ)$-computable) and for all rational numbers $a$ with $0<a<x$, ($f(a)$ exists and) $a<f(a)<x$.

\begin{lem}\label{l3}
There is a positive real number that has no lower separation function.
\end{lem}

\pproof
  We define such a number $z$ by brute force diagonalization.
Let $f_1,f_2,\ldots$ be a sequence of all computable partial functions $f\pf\IQ\to\IQ$.
Let $(a_0;b_0):=(0;1)$ and for $i=1,2,\ldots$ define intervals $(a_i;b_i)$  as follows.
Find some rational number $a$ such that $a_{i-1}<a<f_i(a)<b_{i-1}$ and define
$(a_i;b_i):=(a,f_i(a))$, if no such $a$ exists define $(a_i;b_i):=(a_{i-1};b_{i-1})$.
There is some positive real number $z\in\bigcap_i(a_i;b_i)$. Suppose, $f_k$ is a lower separating function for $z$. Since $z\in (a_{k-1};b_{k-1})$ there is some $a$ such that
$a_{k-1}<a<f_k(a)<z<b_{k-1}$ and $a_k, b_k$ are chosen such that $a_{k-1}<a_k<f_k(a_k)<b_{k-1}$ and $b_k=f_k(a_k)$. If $a_k<z$ then $b_k=f_k(a_k)<z$,
hence $z\not\in(a_k;b_k)$, if $z\leq a_k$ then $z\not\in(a_k;b_k)$ as well.
But by assumption $z\in\bigcap_i(a_i;b_i)$. Therefore, $f_k$ cannot be a lower separating function for $z$.
\qq
\medskip 

\begin{exa}\label{e3}\rm ($T_4$ and $SCT_2$ but not $WCT_3$)
Let ${\bf R}=(\IR,\tau,\beta,\nu)$ be the computable real line from Example~\ref{e8}(\ref{e8a}).
For $c\in\IR$ define ${\bf R}_c=(\IR,\tau,\beta_c,\nu_c)$ by $\nu_c(0w):=\nu(w)$,
and $\nu_c(1w):=\nu(w)\cap(-\infty;c)$. Then ${\bf R}_c$ is a computable topological space.
Let $\delta_c$ be the (canonical or inner) representation of $\IR$ for ${\bf R}_c$
\cite[Definition~5.1]{WG09}.
Since ${\bf R}_c$ has the same topology as ${\bf R}$ it is $T_i$ for $i=0,\ldots,4$.
The computable real line ${\bf R}$ is $SCT_2$. Let $H$ satisfy (\ref{f69}) and (\ref{f70})
for ${\bf R}$.
Then $H_c:=\{(0v,0w)\mid (v,w)\in H\}$ satisfies (\ref{f69}) and (\ref{f70}) for ${\bf R}_c$.
Therefore, ${\bf R}_c$ is $SCT_2$.

 Let $c>0$ be a real number that has no
lower separation function (Lemma~\ref{l3}). Suppose ${\bf R}_c$ is $WCT_3$.
Let $t^w_3$ be the computable function from Definition~\ref{d4} for ${\bf R}_c$. Then:\\
-- the function $a\mapsto a$ for rational $0<a<c$ is $(\nu_\IQ,\delta_c)$-computable,\\
-- $\nu_c(w_0)=(0;c)$ for some $w_0\in\s$,\\
-- for every $x\in(0;c)$, $t^w_3$ maps$(x,(0;c))$ to some $U\in\beta_c$ such that $x\in U\In \overline U\In (0;c)$.\\
-- $U\mapsto\sup U$ for $U\in\beta_c$ such that $\overline U\In (0;c)$ is $(\nu_c,\nu_\IQ)$-computable.

There is a $(\nu_\IQ,\nu_\IQ)$-computable multi-function $h$ mapping each rational
number $0<a<c$ to some rational number $b$ such that $a<b<c$\ :
From a $\nu_\IQ$-name of $a$ compute a $\delta_c$-name of $a$. By $t^w_3$, from $a$ and $(0;c)$ compute some $U\in\beta_c$ such that $a\in U\In\overline U\In (0;c)$. From $U$ compute $b:=\sup U\in\IQ$. Then $a<b<c$.
Since there is an injective notation equivalent to $\nu_\IQ$, the function $h$ is single-valued. Therefore, $h$ is a
lower separation function for $c$. Contradiction.
\qq
\end{exa}

The above space is $T_4$, $T_3$ and $SCT_2$ but not $WCT_3$, $CT_3$ or $SCT_3$.
Finally we separate $\rm CTy$ from $\rm SCT_3$. (The example in
\cite{Wei09} in \cite{BHK09} for separating $\rm CT_3$ from $\rm SCT_3$ is not correct.)

\begin{exa}\label{e14}\rm ($\rm D$ and $\rm CTy$ but not $\rm SCT_3$)
Let $X:=\IN$ and let $A\In \IN$  be the set defined below.
Let $\tau$ be the discrete topology on $\IN$ and define  a notation $\nu$ of a base $\beta$ of $\tau$ by $\nu(12^j):=\{j\}$, $\nu(2):=A$ and $\nu(02^j):=\{j\}\cap A$.
Then ${\bf X}:=(\IN,\tau,\beta,\nu)$ is a computable topological space. Let $\delta$ be the canonical representation of the points of $\bf X$ \cite{WG09}.

We show that $\bf X$ is $CT\hspace{-.3ex}y'$. \\
Let $p_0,p_1\in\om$ be computable sequences such that $\rho(p_0)=0\in\IR$ and $\rho(p_1)=1\in\IR$. There is a machine $M$ that in input $(p,w,q)\in\om\times\s\times\om$ searches for $i,j\in\IN$ such that $12^i\ll p$ and
$12^j\ll q$ and then writes $p_0$ if $i=j$ and $p_1$ else. Then for all $p,q\in\dom(\delta)$ and all $w\in\s$, $\rho\circ f_M(p,w,q)=$ ($0$ if $\delta(p)=\delta(q)$ and $1$ else). Therefore, $f_M$ realizes the function $f:X\times \beta\times X\to\IR$ such that $f(x,W,y)=$  ($0$ if $x=y$ and $1$ else).
By type conversion
\cite[Theorem~3.3.15]{Wei00} the function $(x,W)\mapsto f$ such that $f(y)=$ ($0$ if $x=y$ and $1$ else) is $(\delta,\nu,[\delta\to\rho])$-computable. Furthermore, from
$p\in\dom(\delta)$ we can compute the (unique) $i$ such that $12^i\ll p$. Then
$\delta(p)=i$ and $\{i\}=\nu(12^i)=:U$. Therefore, from $x$ and $W$ such that $x\in W$
we can find $U$ and $f$ such that the conditions for $\rm CTy'$ hold true.

We define $A\In \IN$.
Let $K\In\IN$ be a set with non r.e. complement. Let $A\cap (2\IN+1):=2K+1$.
Define $A\cap 2\IN$ as follows. Let $\gamma_i$ be the $i$th computable function
$f\pf \s\times\s\to\s$  ($i=0,1,2,\ldots$) and let $\lambda_i$ be the $i$th computable function $f\pf \s\times\s\to\om$. For $n=\langle i,k\rangle$ ($\langle\;\rangle$ is the Cantor pairing function) define the position of $2n$ as
follows by diagonalization.
\[\begin{array}{rlll}
\mbox{if}&(12^{2n},2)\not\in\dom(\gamma_i) \mbox{ and } (02^{2n},2)\not\in\dom(\gamma_i)&
\mbox{then }& 2n\in A,\\
\mbox{if}&(12^{2n},2)\in\dom(\gamma_i)  &
\mbox{then }& 2n\not\in A,\\
\mbox{if}&(12^{2n},2)\not\in\dom(\gamma_i) \mbox{ and } (02^{2n},2)\in\dom(\gamma_i)\\
\mbox{then} \\
&\mbox{if}\ \ 02^{2n}\not\ll\lambda_k(02^{2n},2)\ \mbox{and}\ 12^{2n}\not\ll\lambda_k(02^{2n},2)&\mbox{then}& 2n\not \in A,\\
&&\mbox{else} &2n\in A.
\end{array}\]
Suppose, $\bf X$ is $SCT_3$. Let $R$ be the r.e. set and let $r$ be the computable function such that (\ref{f39}) and (\ref{f40}) hold true. Then there are $i,k\in\IN$ such that $R=\dom(\gamma_i)$ and $r=\lambda_k$.

Suppose, $(2,2)\in\dom(\gamma_i)$. Then by (\ref{f40}), $A=\nu(2)\In\psi^-\circ\lambda_k(2,2)\In\nu(2)$,
hence $\IN\setminus A=\bigcup\{\nu(v)\mid v\in V\}$ for an r.e. set $V\In\dom(\nu)$.
Since $2\not\in V$ and $\nu(02^l)=\emptyset$ for $l\not \in A$,
$m\not\in K\iff 2m+1\not\in A\iff  12^{2m+1}\in V$, hence the complement of $K$ is r.e.
(contradiction).
Therefore,  $(2,2)\not\in\dom(\gamma_i)$.

We show that $\gamma_i$ and $\lambda_k$ cannot operate correctly for $n:=\langle i,k\rangle$.
\medskip

\noindent{\bf Case} $(12^{2n},2)\not\in\dom(\gamma_i)$ and $(02^{2n},2)\not\in\dom(\gamma_i)$:
Since $(2,2)\not\in\dom(\gamma_i)$, $2n\not\in\nu(2)$ by (\ref{f39}). But $2n\in A=\nu(2)$ by the definition of $A$. Contradiction.

\medskip
\noindent{\bf Case} $(12^{2n},2)\in\dom(\gamma_i)$: Then $2n\in A=\nu(2)$ by (\ref{f39}). But $2n\not \in A=\nu(2)$ by the definition of $A$. Contradiction.

\medskip
\noindent{\bf Case} $(12^{2n},2)\not\in\dom(\gamma_i)$ and $(02^{2n},2)\in\dom(\gamma_i)$:  By (\ref{f39}),
\begin{eqnarray}\label{f1} \{2n\}\cap A=\nu(02^{2n})\In \psi^-\circ \lambda_k(02^{2n},2)\In \nu(2)=A.
\end{eqnarray}
\indent Suppose $2n\in A$. Then $2n \in\psi^-\circ \lambda_k(02^{2n},2)\In A$,
hence $2n\not\in \theta\circ \lambda_k(02^{2n},2)$.
Therefore,
$02^{2n}\not\ll\lambda_k(02^{2n},2)$ and $12^{2n}\not\ll\lambda_k(02^{2n},2)$.
Then $2n\not\in A$ by the definition of $A$. Contradiction.

 Suppose $2n\not\in A$. By (\ref{f1}), $2n\in \theta\circ \lambda_k(02^{2n},2)$, hence
$02^{2n}\ll\lambda_k(02^{2n},2)$ or $12^{2n}\ll\lambda_k(02^{2n},2)$. Then $2n\in A$ by the definition of $A$. Contradiction.

Therefore, the space $\bf X$ is not $\rm SCT_3$.\qq
\end{exa}

We summarize the counterexamples.

\begin{thm}\label{t1} The following implications are false:
\begin{eqnarray}
\label{t1a} CT_2+D  \  \Longrightarrow  & SCT_2 & \mbox{(\cite[Example~5]{Wei10})},\\
\label{t1b} CT_0+CT_4  \ \Longrightarrow &  T_1 & \mbox{(Example~\ref{e8}(\ref{e8b}))},\\
\label{t1c} SCT_2  \   \Longrightarrow  & T_3 & \mbox{(Example~\ref{e4})},\\
\label{t1d} SCT_2+T_4  \   \Longrightarrow &  WCT_3 & \mbox{(Example~\ref{e3})},\\
\label{t1e} CTy+D   \   \Longrightarrow  & SCT_3  & \mbox{(Example~\ref{e14})}.
\end{eqnarray}
\end{thm}

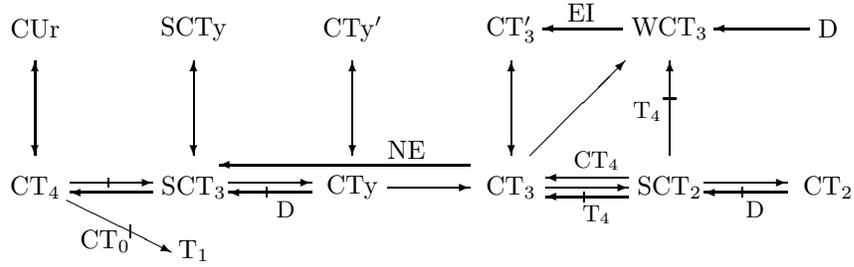
\begin{figure}[htbp]
\setlength{\unitlength}{1.2pt}
\linethickness{0.4pt}
\begin{picture}(180,79)(-25,-20)

\put(-50,0){\makebox(0,0)[cc]{\small $\rm CT_4$}}
\put(-39,1.5){\vector(1,0){26}}
\put(-13,-1.5){\vector(-1,0){26}}
\put(-27,0){\line (0,1){3}}

\put(0,0){\makebox(0,0)[cc]{\small $\rm SCT_3$}}
\put(11,1.5){\vector(1,0){26}}
\put(37,-1.5){\vector(-1,0){26}}
\put(23,-3){\line (0,1){3}}
\put(29,-7){\makebox(0,0)[cc]{\footnotesize $\rm D$}}

\put(50,0){\makebox(0,0)[cc]{\small $\rm CTy$}}
\put(61,0){\vector(1,0){26}}
%\put(87,-1.5){\vector(-1,0){26}}
%\put(73,-3){\line (0,1){3}}
%\put(77,-7){\makebox(0,0)[cc]{\footnotesize $\rm D$}}

\put(100,0){\makebox(0,0)[cc]{\small $\rm CT_3$}}
\put(111,0){\vector(1,0){26}}
\put(137,-3){\vector(-1,0){26}}
\put(137,3){\vector(-1,0){26}}
\put(123,-4.5){\line (0,1){3}}
\put(127,-8.5){\makebox(0,0)[cc]{\footnotesize $\rm T_4$}}
\put(127,8.5){\makebox(0,0)[cc]{\footnotesize $\rm CT_4$}}

\put(150,0){\makebox(0,0)[cc]{\small $\rm SCT_2$}}
\put(161,1.5){\vector(1,0){26}}
\put(187,-1.5){\vector(-1,0){26}}
\put(173,-3){\line (0,1){3}}
\put(177,-7){\makebox(0,0)[cc]{\footnotesize $\rm D$}}

\put(200,0){\makebox(0,0)[cc]{\small $\rm CT_2$}}

%obere Zeile
\put(-50,50){\makebox(0,0)[cc]{\small $\rm CUr$}}
\put(0,50){\makebox(0,0)[cc]{\small $\rm SCTy$}}
\put(50,50){\makebox(0,0)[cc]{\small $\rm CTy'$}}
\put(100,50){\makebox(0,0)[cc]{\small $\rm CT_3'$}}
\put(150,50){\makebox(0,0)[cc]{\small $\rm WCT_3$}}
\put(200,50){\makebox(0,0)[cc]{\small $\rm D$}}

\put(194,50){\vector(-1,0){30}}

\put(135,50){\vector(-1,0){25}}
\put(122,55){\makebox(0,0)[cc]{\small $\rm EI$}}

\put(87,7){\vector(-1,0){79}}
\put(67,12){\makebox(0,0)[cc]{\small $\rm NE$}}

\put(106,10){\vector(1,1){30}}

\put(150,10){\vector(0,1){30}}
\put(148,28){\line (1,0){4}}
\put(143,24){\makebox(0,0)[cc]{\footnotesize $\rm T_4$}}

\put(100,10){\vector(0,1){30}}
\put(100,40){\vector(0,-1){30}}

\put(50,10){\vector(0,1){30}}
\put(50,40){\vector(0,-1){30}}

\put(0,10){\vector(0,1){30}}
\put(0,40){\vector(0,-1){30}}

\put(-50,10){\vector(0,1){30}}
\put(-50,40){\vector(0,-1){30}}
%unten

\put(-40,-4){\vector(2,-1){33}}
\put(0,-20){\makebox(0,0)[cc]{\small $\rm T_1$}}
\put(-20,-16){\line (0,1){4}}
\put(-28,-17){\makebox(0,0)[cc]{\small $\rm CT_0$}}

\end{picture}
\caption{Logical relations between the computable separation axioms.} \label{fig3}
\end{figure}

Further false implications can be obtained by transitivity of \, "$\Longrightarrow$", for example, $\rm CT_4\Longrightarrow SCT_3$ is false
by (\ref{t1b}) since $\rm SCT_3\Longrightarrow T_1$.
Figure~\ref{fig3} shows the positive and negative results that we have proved.
``$A\longrightarrow B$'' means $A\Longrightarrow B$,
``$A\longrightarrow$~\hspace{-3.5ex}\raisebox{1ex}{\footnotesize $C$}\hspace{1.ex} $B$'' means $A\wedge C \Longrightarrow B$,
``$A\not\longrightarrow B$'' means that we have constructed  a
 computable topological space for which $A\wedge \neg B$, and
``$A\not\longrightarrow$~\hspace{-3ex}\raisebox{1ex}{\footnotesize $C$}\hspace{1.ex} $B$'' means that we have constructed a computable topological space for which $(A\wedge C)\wedge \neg B$. ${\rm EI}$ abbreviates ``$U\cap V =\emptyset$ is $(\nu,\nu)$-r.e.'' and ${\rm NE}$ abbreviates ``$U\neq\emptyset$ is $\nu$-r.e.''.

A number of implications have not yet been proved or disproved, for example,
 \[\rm \begin{array}{ccl}
WCT_3 &  \Longrightarrow  & CT_3 \ (SCT_2, CT_2,CT_0,WCT_0) ,\\
CT_3 &  \Longrightarrow  & CTy \ (CT_4),\\
CT_4+T_3 &  \Longrightarrow  & (SCT_3, CT_3,SCT_2,CT_2,CT_0,WCT_0),\\
CT_4+ SCT_2 &\Longrightarrow  &SCT_3 \ (CTy).
\end{array}\]
(The axioms $\rm CT_0$ and $\rm WCT_0$ are defined in \cite{Wei10}.)
A difficulty arises from the fact that for $T_1$-spaces (where the singleton sets are closed) the function $x\mapsto\{x\}$ is $(\delta,\psi^+)$-computable but in general not $(\delta,\psi^-)$-computable. In our computable separation axioms, however, we use only the outer representation $\psi^-$ for the closed sets.

\section{Computable Metrization}\label{secs}

For a metric space $(X,d)$ the open balls with rational radius and center from a dense set are a basis of a topology, the topology generated by it \cite{Eng89}.
A topological space $(X,\tau)$ is metrizable, iff it is generated by a metric space.
A second countable space is metrizable iff it is $T_3$ (remember $\rm T_3\Longrightarrow T_2$) \cite[Theorem~4.2.9]{Eng89}. Every second countable metrizable space has an at most countable dense subset $A$. Therefore, it can be enriched by a notation $\alpha$ of this dense set.
The following definitions are essentially from \cite{Wei93,Wei00}.

\begin{defi}\label{d1}\hfill
\begin{enumerate}[(1)]
\item\label{d1c} An effective metric space is a tuple ${\bf M}=(X,d,A,\alpha)$ such that $(M,d)$ is a metric space and $\alpha$ is a notation of a set $A\In X$ which is dense in $X$.

\item\label{d1e} The Cauchy representation $\delta_C$ of an effective metric space $\bf M$ is defined by $\delta_C(p)=x$ iff there are words $u_0,u_1,\ldots\in\dom(\alpha)$ such that $p=\iota(u_0)\iota(u_1)\ldots$ and $d(x,\alpha(u_i))\leq 2^{-i}$ for all $i\in\IN$.

\item\label{d1d} The effective topological space \cite[Definition~4]{WG09} associated with the effective metric space is the tuple ${\bf X}=(X,\tau,\beta,\nu)$ such that $\nu$,  $\nu\langle u,s\rangle:=B(\alpha(u),\nu_\IQ(s))$, is the
canonical notation of the set $\beta$ of all open balls with center from $A$ and rational radius and $\tau$ is the smallest topology containing~$\beta$.

\item\label{d1a}  An upper semi-computable (lower semi-computable) metric space is an effective metric space such that $\dom(\alpha)$ is recursive and $ d(a,b)<s$
($ s<d(a,b)$) is
$(\alpha,\alpha,\nu_\IQ)$-r.e.
\item\label{d1b} A computable metric space is an effective metric space such that
$\dom(\alpha)$ is recursive and  $ r<d(a,b)<s$ is $(\nu_\IQ,\alpha,\alpha,\nu_\IQ)$-r.e.
\end{enumerate}
\end{defi}

\noindent Notice that $ d(a,b)<s$ is $(\alpha,\alpha,\nu_\IQ)$-r.e. iff the distance on $A$ is $(\alpha,\alpha,\rho_>)$\bb computable iff $d$ is $(\delta_C,\delta_C,\rho_>)$-computable \cite[Example~1]{WG09}\cite{Wei00}, and $r< d(a,b)<s$ is $(\nu_\IQ,\alpha,\alpha,\nu_\IQ)$-r.e. iff the distance on $A$ is $(\alpha,\alpha,\rho)$-computable iff the distance is $(\delta_C,\delta_C,\rho)$-computable (Example~\ref{e8}(\ref{e8a})) \cite{Wei00}.
Since for every notation $\alpha\pf \s\to A$ with r.e. domain there is a notation $\alpha'\pf\s\to A$ with recursive domain such that $\alpha\equiv\alpha'$, allowing r.e. domains in Definition~\ref{d1}(\ref{d1a}) and~(\ref{d1b}) is no proper generalization.

\begin{thm}\label{t3}For every effective metric space ${\bf M}=(X,d,A,\alpha)$ with Cauchy representation $\delta_C$ and its associated effective topological space ${\bf X}=(X,\tau,\beta,\nu)$ with canonical representation $\delta$,
\begin{enumerate}[\em(1)]
\item\label{t3a}if $\;\bf M$ is upper semi-computable, then ${\bf X}$ is a computable topological space.

\item\label{t3b}$\delta\leq \delta_C$; \ $\delta_C\leq_t \delta$; \ $\delta_C\leq \delta$
if $\;\bf M$ is upper semi-computable,
\item\label{t3c}   if  $\bf M$ is a computable metric space, then all the separation axioms from Definition~\ref{d4} hold true for ${\bf X}$.
\end{enumerate}
\end{thm}

\noindent The first two items of this theorem  differ slightly from \cite[Theorem~8.1.4]{Wei00} since {\em computable topological space} is defined differently.
\medskip

\pproof (\ref{t3a}) We must show that intersection is $(\nu,\nu,\theta)$-computable. Observe that\\
$B(a_1,r_1)\cap B(a_2,r_2)=\bigcup\{B(a,r)\mid d(a_1,a)<r_1-r \wedge d(a_2,a)<r_2-r\}$.
Since $ d(a,b)<s$ is r.e., there is an r.e. set $S$ such that
$\nu(u)\cap\nu(v)=\bigcup \{\nu(w)\mid (u,v,w)\in S\}$.

(\ref{t3b}) If $\delta(p)=x$, then $p$ is a list of all $\langle u,v\rangle$ such that
$d(x,\alpha(u))<\nu_\IQ(v)$. Therefore, there is a machine $M$ which from $p\in\dom(\delta)$ computes a sequence $\iota(u_0)\iota(u_1)\ldots$ such that for all $i$, $d(x,\alpha(u_i))\leq2^{-i}$.

 Now let $\bf M$ be upper semi-computable.
If $\delta_C(p)=x$ then $p=\iota(u_0)\iota(u_1)\ldots$ such that for all $i$, $d(x,\alpha(u_i))\leq2^{-i}$. Observe that
$x\in B(a,r)\iff (\exists i)\, d(a,\alpha(u_i))<r-2^{-i}$. Since $ d(a,b)<s$ is r.e., from $p$ we can compute a list of all $w$ such that $x\in\nu(w)$.

If $\bf M$ is not upper semi-computable then there are ``oracles'' $q,q'\in\om$ such that
$\dom(\alpha)$ is recursive in $q$ and
$d(a,b)<s$ is r.e. in $q'$. Using the oracles $q,q'$, there is a machine translating
$\delta_C$ to $\delta$. The function $f_M$ computed by this machine is continuous \cite{Wei00}.

(\ref{t3c}) By Theorem~\ref{t6} it suffices to prove $\rm SCT_3$.
Let $R:=\{(\langle u,v\rangle,\langle u',v'\rangle)\mid
\langle u,v\rangle,\langle u',v'\rangle\in\dom(\nu), \ \
d(\alpha(u),\alpha(u'))+\nu_\IQ(v)<\nu_\IQ(v')\}$.
Then (\ref{f39}) and $\overline{\nu(\langle u,v\rangle)}\In \nu(\langle u',v'\rangle)$\\
for $(\langle u,v\rangle,\langle u',v'\rangle)\in R$. We compute a $\psi^-$-name of this closure.
There is a machine that on input
$(\langle u,v\rangle,\langle u',v'\rangle)\in R$ lists all $\langle w,w'\rangle\in\dom(\nu)$ such that $d(\alpha(u),\alpha(w))>
\nu_\IQ(v) +\nu_\IQ(w')$. Then (\ref{f40}) holds true for the function $r:=f_M$.
\qq\\

Since for an effective topological space $(X,\tau,\beta,\nu)$, $\delta$ is an admissible representation, by Theorem~\ref{t3}(\ref{t3b}) the Cauchy representation is admissible, that is, it is continuous and $\delta\leq_t\delta_C$ for every continuous representation of~$X$ \cite{Wei00}.

In general, we are interested in metric spaces $(X,d)$ with representation $\delta\pf\om\to X$ such that the distance is at least $(\delta,\delta,\rho_>)$-continuous. In this case
the metric space is separable and the representation $\delta\pf\om\to X$ is continuous \cite[Lemma~8.1.1]{Wei00}. By adding a notation of a dense set we obtain an effective metric space with Cauchy representation $\delta_C$.
Then $\delta\leq_t\delta_C$, since the Cauchy representation is admissible.

We call a metric on a computable topological space ${\bf X}=(x,\tau,\beta,\nu)$ with canonical representation $\delta$ of the points lower semi-computable, if it is $(\delta,\delta,\rho_<)$-computable and computable, if it is $(\delta,\delta,\rho)$-computable.

\begin{thm}\label{t5} Let $\bf X$ be a computable topological space.
\begin{enumerate}[\em(1)]
\item\label{t5a} Suppose some lower semi-computable metric $d$ generates the topology of $\bf X$. Then $\bf X$ is $SCT_2$.
\item\label{t5b}\cite{Sch98,GSW07} Suppose $\bf X$ is $SCT_3$. Then its topology is generated by some computable metric.
\end{enumerate}
\end{thm}

Theorem~\ref{t5}(\ref{t5b}) has been proved in \cite{Sch98}.
The shorter proof in \cite{GSW07} assumes  $U\neq\emptyset$ for $U\in\beta$ but actually does not need this condition. We include a proof, since parts of it will be used in the proof of the next theorem.\\

\pproof

(\ref{t5a})  By \cite[Theorem~11]{WG09} there is a computable function $g$ such that $\nui(w)=\theta\circ g(w)$. There is a machine $M$ such that $f_M$ realizes the distance function w.r.t. $(\delta,\delta,\rho_<)$.
Let $H$ be the set of all $(u,v)$ for which there are $v_1,v_2\in\dom(\nufs)$ and $v_3,v_4$ such that\\
(a) the machine $M$ on input $(v_11^\omega,v_21^\omega)$ writes in at most $\max(|v_1|,|v_2|)$ steps the word $v_3$, $v_4\ll v_3$  and $\nu_\IQ(v_4)>0$, and\\
(b) $u\ll g(v_1)$ and $w\ll g(v_2)$.\\
The set $H$ is r.e. We prove (\ref{f69}) amd (\ref{f70}).

Suppose $\delta(p)=x\neq y=\delta(q)$. Since $d(x,y)>0$ there are $v_1\pref p$, $v_2\pref q$, $v_3$ and $v_4$ such that the machine $M$ on input $(v_11^\omega, v_21^\omega)$ writes in at most $\max(|v_1|,|v_2|)$ steps the word $v_3$, $v_4\ll v_3$ and $\nu_\IQ(v_4)>0$.
Since $x\in \bigcap \nufs(v_1)$ and $y\in \bigcap \nufs(v_2)$, there are $u\ll g(v_1)$ and $w\ll g(v_2)$ such that $x\in\nu(u)$ and $y\in \nu(w)$. By the definition of $H$, $(u,w)\in H$. This proves  (\ref{f69}).

Suppose $(u,w)\in H$. If $\nu(u)=\emptyset$ or $\nu(w)=\emptyset$, then $\nu(u)\cap\nu(w)=\emptyset$. Suppose $\nu(u)\neq\emptyset$ and $\nu(w)\neq\emptyset$.
Then there are words $v_1,v_2,v_3$ and $v_4$ such that,
the machine $M$ on input $(v_11^\omega,v_21^\omega)$ writes in at most $\max(|v_1|,|v_2|)$ steps the word $v_3$, $v_4\ll v_3$,  $\nu_\IQ(v_4)>0$, $u\ll g(v_1)$ and $w\in g(v_2)$.
Since $\nu(u)\In \bigcap \nufs(v_1)$ and $\nu(w)\In \bigcap \nufs(v_2)$, for every $x\in\nu(u)$ and every $y\in\nu(w)$ there are sequences $p',q'$ such that
$x=\delta(v_1p')$ and $y=\delta(v_2q')$. Then the machine $M$ on input $(v_1p',v_2q')$ writes in at most $\max(|v_1|,|v_2|)$ steps the word $v_3$ such that $v_4\ll v_3$ and $\nu_\IQ(v_4)>0$. Therefore, $d(x,y)>0$ for every $x\in\mu(u)$ and $y\in\nu(w)$.  This proves  (\ref{f70}).

(\ref{t5b})
Since $R=\dom(r)$ is r.e., it has a computable numbering $(u_i,v_i)_{i\in\IN}$.
By Theorem~\ref{t6}, the  Urysohn multi-function~$t_{\rm Ur}$ has a computable
$(\psi^-,\psi^-,[\delta \to \rho])$-realization $h\pf \om\times\om\to\om$.

The function $U\mapsto U^c$ for $U\in\beta$ has a computable $(\nu,\psi^-)$-realization $g\pf \s\to \om$. For $i\in\IN$ define $f_i:X\to\IR$, $d_i:X\times X\to \IR$ and $d_:X\times X\to \IR$ by
\begin{eqnarray}
\label{f2}f_i &:=& [\delta\to\rho]\circ h(r(u_i,v_i),g(v_i))\\
\label{f3}d_i(x,y)&:=&|f_i(x)-f_i(y)|\\
\label{f4}d(x,y)&:=&\sum_i 2^{-i}d_i(x,y)
\end{eqnarray}
Then for every $i$, $f_i$ is a continuous function such that $\range(f_i)\In [0;1]$, $f(x)=0$ for $x\in \nu(u_i)$ and $f(x)=1$ for $x\not\in\nu(v_i)$, and
$d_i$ is a continuous pseudometric  on $(X,\tau)$ bounded by $1$ such that $d_i(x,y)=1$ for $x\in \nu(u_i)$ and $y\not\in\nu(v_i)$.

Let $A$ be closed and non-empty and $x\not\in A$.
Then there is some $i$ such that $A\In (\nu(v_i))^c$ and
$x\in \nu(u_i)\In \nu(v_i)$. Then $d_i(x,A):=\inf\{a\in A\mid d_i(x,a)\}=1$.
By \cite[Lemma~4.4.6]{Eng89}, $d$ is a metric which generates the topology $\tau$.

Since $i\mapsto f_i$ is $(\nu_\IN,[\delta\to\rho])$-computable, the metric $d$ is
$(\delta,\delta,\rho)$-computable.
\qq\\

The condition in the metrization theorem~\ref{t5}(\ref{t5b}) is $SCT_3$. We do not know whether $STy$ or $CT_3$ are sufficient to prove the metrization  theorem.

For a computable metric space a dense set of computable points is needed. In general  a space with computable metric does not have computable points but its metric completion may have computable points (example: the restriction of the computable real line
(Example~\ref{e8}(\ref{e8a})) to the non-computable real numbers). We will show that the metric space constructed in the proof of Theorem~\ref{t5} can be completed to a computable metric space, if $\{u\in\dom(\nu)\mid \nu(u)\neq\emptyset\}$ is r.e..

For a pseudo-metric $d$ and sets $A,B$ we define the diameter and the distance of sets as usual: \  ${\rm dm}(A):=\sup\{d(x,y)\mid x,y\in A\}$, \
$d(A,B):=\inf \{d(x,y)\mid x\in A,\ y\in B\}$. The triangle inequality generalizes to
\begin{eqnarray}\label{f22}d(A,C)\leq d(A,B) +d(B,C)+ {\rm dm}(A)  +{\rm dm}(B) +{\rm dm}(C)\,.
\end{eqnarray}

We define a computable version of {\em homeomorphic embedding} \cite[Section~2.1]{Eng89}.
We will construct a computable metric space such that original computable topological space can be computably embedded into it.

\begin{defi}\label{d3}
For represented spaces $(X,\delta)$ and $(X',\delta')$, a {\em computable embedding} is an injective function $f:X\to X'$ such that $f$ is $(\delta,\delta')$-computable and
$f^{-1}$ is $(\delta',\delta)$-computable.
\end{defi}

For computable topological spaces the standard representations are admissible, hence relatively computable functions are continuous.  In this case a computable embedding is a homeomorphic embedding.

\begin{thm}\label{t8}\cite{GSW07} Let ${\bf X}=(X,\tau,\beta,\nu)$ be a computable topological space such that $\rm CT_3$ is true and the set $\{u\in\dom(\nu)\mid \nu(u)\neq\emptyset\}$ is r.e.
Then there is a computable embedding of $\bf X$ into a computable metric space
${\bf M}=(M,d_M,A,\alpha)$ (where for $\bf X$ we consider the standard representation and for $\bf M$ the Cauchy representation).
\end{thm}

In \cite{GSW07} the theorem has been proved for $SCT_3$ spaces with non-empty base sets.
First we show that the assumptions in Theorem~\ref{t8} are sufficient and then present a
proof that uses ideas from \cite{GSW07} but is more transparent and much simpler.
\bigskip

\pproof
By \cite[Lemma~25]{WG09} there is a computable topological space
${\bf X'}=(X,\tau,\beta',\nu')$ equivalent to $\bf X$ such that $\nu'(u)\neq\emptyset$ for all $u\in\dom(\nu')$. By \cite[Theorem~22]{WG09} equivalent means $\delta\equiv\delta'$, hence the identity is a computable embedding of $\bf X$ into
$\bf X'$. By Lemma~\ref{l5} we may assume w.l.o.g. that $\nu(u)\neq\emptyset$ for all $u\in\dom(\nu)$.

By Theorem~\ref{t7} the space ${\bf X}$ is $SCT_3$. For $i\in\IN$ let $f_i$ be the level function and $d_i$ the pseudo-metric  and let $d$ be the metric defined in the proof of Theorem~\ref{t5} ((\ref{f2}), (\ref{f3}), (\ref{f4}))  with diameters $\dm_i$ and $\dm$, respectively. Remember that $i\mapsto f_i$ is $(\nu_\IN,[\delta\to\rho])$-computable.

In the following we use nested sequences of non-empty open sets instead of Cauchy-sequences of points for completion.\newpage

\begin{prop}\label{p1}\hfill
\begin{enumerate}[\em(1)]
\item \label{p1a}
The multi-function $g_1$ mapping every $(W,i,n)$ such that  $W\in\beta$ and $i,n\in\IN$ to some $(U,a)$ such that $U\in\beta$, $a\in\IQ$, $U\In W$ and $f_i[U]\In (a-2^{-n};a+2^{-n})$ is computable.
\item \label{p1b}
The multi-function $g_2$ mapping every $(x,W,i,n)$ such that $x\in W\in\beta$ and $i,n\in\IN$  to some $(U,a)$ such that $x\in U\in\beta$, $a\in\IQ$, $U\In W$ and $f_i[U]\In (a-2^{-n};a+2^{-n})$ is computable.
\end{enumerate}
\end{prop}

\pproof (Proposition~\ref{p1})

(\ref{p1a})
By the statement $\overrightarrow{\delta_1}\equiv \overrightarrow{\delta_3}$ in
\cite[Theorem~29]{WG09} and (\ref{f2}) and since intersection on open sets is computable \cite[Theorem~11]{WG09}, there is a computable function $g$ mapping $(i,n,a,u')$ ($i,N\in\IN$, $a\in\IQ$, $u'\in\dom(\nu)$) to some $q\in\om$ such that
$f_i^{-1}[(a-2^{-n};a+2^{-n})]\cap\nu(u')=\theta(q)$. There is a machine $M$ that on input $(w,i,n)$ searches for some $a\in\IQ$ and $u,u'\in \dom(\nu)$ such that
$(u',w)\in R$, $\nu(u)\neq\emptyset$ and $u\ll g(i,n,a,u')$ and then writes $(u,a)$.

 There is some  $y\in W=\nu(w)$. Then there is some $u'$ such that $(u',w)\in R$ and $y\in\nu(u')$. There is some $a$ such that $f_i(y)\in (a-2^{-n};a+2^{-n})$. Since
$y\in f_i^{-1}[(a-2^{-n};a+2^{-n})]\cap\nu(u')$ there is some $u$ such that
$u\ll g(i,n,a,u')$. Therefore, the machine $M$ on input $(w,i,n)$ succeeds to write some $(u,a)$. In this case, $\nu(u)\In\nu(u')\In\nu(w)$ and $f_i[\nu(u)]\In (a-2^{-n};a+2^{-n})$. This proves the first statement.

(\ref{p1b})  Let the machine from the above proof search for some $u$ such that
additionally $x\in\nu(u)$.
\hfill$\Box$(Proposition~\ref{p1})
\medskip

We define a computable metric space ${\bf M}=(M,d_M,A,\alpha)$ as the constructive completion of a computable notated pseudometric space ${\bf A}'=(A',d',\alpha')$
\cite[Definition~8.1.5]{Wei00} which will be constructed now.
Let $A':=\dom(\nu)$ and $\alpha'(u):=u$ for $u\in A'$.
Iterating a  computable realization of the multi-function $g_1$ from Proposition~\ref{p1}(\ref{p1a}) for every $w\in A'$ we can compute a sequence
$((u_{wk},a_{wk}))_{k\in\IN}$ (where $(u_{wk},a_{wk})\in A'\times \IQ$) such that
for $k=\langle i,n\rangle$,
\begin{eqnarray}
\label{f5}& \nu(u_{w,k+1})\In \nu(u_{wk})\In \nu(w)\,,\\
\label{f6}&f_i[\nu(u_{wk})]\In (a_{wk}-2^{-n};a_{wk}+2^{-n})\,.
\end{eqnarray}
Then for $v,w\in A'$ define
\begin{eqnarray}\label{f7}
&d'(v,w):=\sup_k d(\nu(u_{vk}),\nu(u_{wk}))\,.
\end{eqnarray}

By (\ref{f5}) the sequence $(\dm\circ \nu(u_{wk}))_{k\in\IN}$ of diameters is decreasing and by (\ref{f6}),
\begin{eqnarray}\label{f9}
\dm_i\circ \nu(u_{w\langle i,n\rangle})\leq 2\cdot 2^{-n}\,.
\end{eqnarray}
Let $x,y\in \nu(u_{w\langle n,n\rangle})$.
Then for all $j\leq n$, $\langle j,n\rangle\leq \langle n,n\rangle$, hence $x,y\in \nu(u_{w\langle j,n\rangle})$ by (\ref{f5}) and therefore,
 $|f_j(x)-f_j(y)|\leq 2\cdot 2^{-n}$ by (\ref{f9}).
 Since ${\rm range}(f_j)\In [0;1]$,
\[\begin{array}{lllll}
d(x,y)&=&\sum_{j\in\IN} 2^{-j}|f_j(x)-f_j(y)|\\
&\leq &\sum_{j\leq n} 2^{-j}|f_j(x)-f_j(y)|+ 2^{-n}\\
&\leq &\sum_{j\leq n} 2^{-j}\cdot2\cdot  2^{-n}+ 2^{-n}\leq 5\cdot  2^{-n}\,,
\end{array}\]
and hence,
\begin{eqnarray}\label{f8}
\hspace{-20ex}\dm(\nu(u_{w\langle n,n\rangle}))&\leq & 5\cdot 2^{-n}\,.
\end{eqnarray}
Since $(\forall k)(\exists n)k\leq \langle n,n\rangle$, the sequence $(\dm\circ \nu(u_{wk}))_{k\in\IN}$ converges to~$0$.\\
For $U_k:=\nu(u_{w_1k})$, $V_k:=\nu(u_{w_2k})$ and $W_k:=\nu(u_{w_3k})$,
\[\begin{array}{lll}
d(U_k,W_k)&\leq& d(U_k,V_k)+d(V_k,W_k)+\dm(U_k)+\dm(V_k)+\dm(W_k)\\
&\leq& d'(w_1,w_2)+d'(w_2,w_3)+\dm(U_k)+\dm(V_k)+\dm(W_k)\,,
\end{array}\]
hence $d'(w_1,w_3)\leq d'(w_1,w_2)+d'(w_2,w_3)$. Therefore, $d'$ is a pseudometric. \\

We will show that $d'$ is computable.
Since we have assumed that the base elements of the space $\bf X$ are not empty, for every $w,k$  there is some $x_{wk}\in \nu(u_{wk})$. Although  we are not able to compute such points we will use their existence.
 For  $m>n$, $x_{v\langle i,m\rangle}\in\nu(u_{v\langle i,m\rangle})\In
\nu(u_{v\langle i,n\rangle})$. Then by (\ref{f6}),
$|a_{v\langle i,m\rangle}- a_{v\langle i,n\rangle}|\leq 2\cdot 2^{-n}$.
Therefore the sequence $(a_{v\langle i,n\rangle})_{n\in\IN}$ converges to some $b_{vi}\in\IR$ such that $|b_{vi}- a_{v\langle i,n\rangle}|\leq 2\cdot 2^{-n}$.
The function $(v,i)\mapsto b_{vi}$ is computable. Furthermore,
$|a_{v\langle i,n\rangle}- a_{w\langle i,n\rangle}| \leq |a_{v\langle i,n\rangle}- b_{vi}|
+ |b_{vi}-b_{wi}| + | b_{wi}- a_{w\langle i,n\rangle}|\leq |b_{vi}-b_{wi}|+ 4\cdot 2^{-n}$
and correspondingly
$ |b_{vi}-b_{wi}| \leq  |a_{v\langle i,n\rangle}- a_{w\langle i,n\rangle}|+4\cdot 2^{-n}$, hence
$|\;|b_{vi}-b_{wi}|- |a_{v\langle i,n\rangle}- a_{w\langle i,n\rangle}|\;|
\leq  4\cdot 2^{-n}$.
By (\ref{f6}) for $k=\langle i,n\rangle$,
$|a_{vk}-a_{wk}|-2\cdot 2^{-n}\leq |f_i(x_{vk})-f_i(x_{wk})|\leq |a_{vk}-a_{wk}|+2\cdot 2^{-n}$, hence  $|\;|a_{vk}-a_{wk}|- d_i(x_{vk},x_{wk})\,|\leq 2\cdot  2^{-n}$.
Therefore,
\begin{eqnarray}\label{f28}|\; |b_{vi}-b_{wi}|- d_i(x_{v\langle i,n\rangle},x_{w\langle i,n\rangle})\,|\leq 6\cdot2^{-n}\,.
\end{eqnarray}
Suppose $m\geq \langle i,n\rangle$. Since $x_{vm}\in\nu(u_{vm})\In \nu(u_{v\langle i,n\rangle})$ and $x_{v\langle i,n\rangle}\in\nu(u_{v\langle i,n\rangle})$, \\
$d_i(x_{v\langle i,n\rangle},x_{vm})\leq 2\cdot 2^{-n}$ by (\ref{f9}) and correspondingly
$d_i(x_{w\langle i,n\rangle},x_{wm})\leq 2\cdot 2^{-n}$, hence
\[|\,d_i(x_{vm},x_{wm})- d_i(x_{v\langle i,n\rangle},x_{w\langle i,n\rangle})\,|\leq  4\cdot 2^{-n}\,,\]
and with (\ref{f28}),
\begin{eqnarray}\label{30}
|\;|b_{vi}-b_{wi}| - d_i(x_{vm},x_{wm})\,|&\leq & 10\cdot2^{-n}\,.
\end{eqnarray}
Then for $N\in\IN$  and $m>\langle N+1,N+6\rangle$,
\begin{eqnarray*}
&&\left |d(x_{vm},x_{wm})-\sum_{i=0}^{N+1}2^{-i}|b_{vi}-b_{wi}|\right |\\
&=& \left |\sum_{i\in\IN}2^{-i}\cdot d_i(x_{vm},x_{wm})-\sum_{i=0}^{N+1}2^{-i}|b_{vi}-b_{wi}|\right |\\
&\leq& \left |\sum_{i=0}^{N+1}2^{-i}\cdot d_i(x_{vm},x_{wm})-\sum_{i=0}^{N+1}2^{-i}|b_{vi}-b_{wi}|\right | +2^{-N-1}\\
&=& \left |\sum_{i=0}^{N+1}2^{-i}\cdot \left (d_i(x_{vm},x_{wm})-|b_{vi}-b_{wi}|\right )\right | +2^{-N-1}\\
&\leq& \sum_{i=0}^{N+1}2^{-i}\cdot \Big |d_i(x_{vm},x_{wm})-|b_{vi}-b_{wi}|\;\Big | +2^{-N-1}\\
&\leq& \sum_{i=0}^{N+1}2^{-i}\cdot 10\cdot 2^{-N-6} +2^{-N-1}\\
&\leq& 20\cdot 2^{-N-6} +2^{-N-1} <2^{-N}\,.\\
\end{eqnarray*}
Since
$d(\nu(u_{vm}),\nu(u_{wm}))\leq d(x_{vm},x_{wm}) \leq d(\nu(u_{vm}),\nu(u_{wm}))+
{\rm dm} (\nu(u_{vm}) + {\rm dm} (\nu(u_{wm})$,
\[\lim_{m\to\infty} d(\nu(u_{vm}),\nu(u_{wm}))=\lim_{m\to\infty} d(x_{vm},x_{wm})=d'(v,w)\]
Therefore by the above estimation,
\[\left |[d'(v,w) -\sum_{i=0}^{N+1}2^{-i}|b_{vi}-b_{wi}|\right| \leq 2^{-N}\]
for all $N$. Since the function $(v,w,N)\mapsto \sum_{i=0}^{N+1}2^{-i}|b_{vi}-b_{wi}|$ is computable, the pseudometric $d'$ on the pseudometric space ${\bf A}'=(A',d',\alpha ')$ is $(\alpha',\alpha',\rho)$-computable.

\smallskip
Let ${\bf M}=(M,d_M,A,\alpha)$ be the constructive completion of  the computable notated pseudometric space ${\bf A}'=(A',d',\alpha')$, see \cite[Definition~8.1.5]{Wei00}. We summarize its definition.
Define a set $S$, a function $d_S:S\times S\to\IR$ and a binary  relation $\sim$ on $S$ as follows:
\begin{eqnarray}
\label{f31} & S  :=  \{(w_0,w_1,\ldots)\mid w_i\in A', \ d'(w_i,w_j)\leq 2^{-i}\ \mbox{for}\ j>i\}\,,\\
\label{f32}&
 d_S((v_0,v_1,\ldots),(w_0,w_1,\ldots)) :=  \lim_{i\to\infty}d'(v_i,w_i)\,,\\
\label{f33}  &(v_0,v_1,\ldots)\sim (w_0,w_1,\ldots) \iff   d_S((v_0,v_1,\ldots),(w_0,w_1,\ldots))=0\,.
\end{eqnarray}
Then define $M:=S\msim$, $d_M:=d_S\msim$, $\alpha(w):=(w,w,w,\ldots)\msim$ for $w\in\dom(\alpha):=\dom(\alpha')=\dom(\nu)=A'$ and $A:=\range(\alpha)$.

The Cauchy representation $\delta_C$ for $\bf M$ is defined by:
$p\in\dom(\delta_C)$ iff there are words $w_0,w_1,\ldots\in\dom(\alpha)$ such that
$p=\iota(w_0)\iota(w_1)\ldots$  and $d'(w_i,w_j)\leq 2^{-i}$ for $j>i$,
and $\delta_C(p)=(w_0,w_1,\ldots)\msim$.
\smallskip

We will define a function $f:X\to M$ and prove that $f$ is well-defined, injective and $(\delta,\delta_M)$-computable and that the partial function $f^{-1}$ is $(\delta_M,\delta)$-computable.
For every $w\in\dom(\nu)$ let $((u_{wk},a_{wk}))_{k\in\IN}$ be the sequence satisfying (\ref{f5}, \ref{f6})  for $k=\langle i,n\rangle$ that has been used for defining the pseudometric space ${\bf A}'$.
\smallskip

Let $\delta(p)=x$. There is a machine $N$ that on input $p\in\om$ first finds some $w\ll p$.
Using a computable realization of the multi-function $g_2$ from Proposition~\ref{p1}(\ref{p1b}) from $p$ and $w$ it computes a sequence
$((v_{pwk},c_{pwk}))_{k\in\IN}$ (where $(v_{pwk},c_{pwk})\in A'\times \IQ$) such that
for $k=\langle i,n\rangle$,
\begin{eqnarray}
\label{f17}&x\in \nu(v_{pw,k+1})\In \nu(v_{pwk})\In \nu(w)\,,\\
\label{f18}&f_i[\nu(v_{pwk})]\In (c_{pwk}-2^{-n};c_{pwk}+2^{-n})\,.
\end{eqnarray}
(compare with (\ref{f5}, \ref{f6})) and writes the sequence
$q:=\iota(v_0)\iota(v_1)\ldots$ where $v_n:=v_{pw\langle n+3,n+3\rangle}$.
In the same way as above from (\ref{f5}, \ref{f6}) from
(\ref{f17}, \ref{f18}) we can conclude
$\dm(\nu(v_{pw\langle n,n\rangle}))\leq 5\cdot 2^{-n}$.
Then
\begin{eqnarray}\label{f19}
x\in \nu(v_{n+1})\In\nu(v_n)&\mbox{and} & \dm(v_n)<2^{-n}\,.
\end{eqnarray}
We show $q\in\dom(\delta_C)$. Suppose $i<j$ and let $n\in\IN$.  Since $\nu(v_j)\In \nu(v_i)$ by (\ref{f5}),
$\nu(u_{v_i\langle n,n \rangle})\In \nu(v_i)$ and
$\nu(u_{v_j\langle n,n \rangle})\In \nu(v_i)$,
hence $d(\nu(u_{v_i\langle n,n \rangle}),
\nu(u_{v_j\langle n,n \rangle}))\leq \dm\circ\nu(v_i)$.
Therefore by (\ref{f7}) and (\ref{f19}),
\[d'(v_i,v_j)=\sup_n d(\nu(u_{v_i\langle n,n \rangle}),
\nu(u_{v_j\langle n,n \rangle}))\leq \dm\circ \nu(v_i)\leq 2^{-i}\,.\]
Therefore, $q\in\dom(\delta_C)$.

Let $\delta(p)=x$, $\delta(p')=x'$, $f_N(p)=q=\iota(v_0)\iota(v_1)\ldots$ and
$f_N(p')=q'=\iota(v_0')\iota(v_1')\ldots$.
By the definition of $d_M$,  $d_M(\delta_C(q),\delta_C(q'))=d_S((v_0,v_1,\ldots),(v_0',v_1',\ldots))
=\lim_{i\to\infty} d'(v_i,v_i')$.

For all $i\in\IN$, $x\in\nu(v_i)$ and $x'\in\nu(v_i')$ by (\ref{f19}) and for all $n\in\IN$, $\nu(u_{v_i\langle n,n  \rangle})\In\nu(v_i)$ and
$\nu(u_{v_i'\langle n,n \rangle})\In\nu(v_i')$ by (\ref{f5}).
For $y\in\nu(v_i)$ and $y'\in\nu(v_i')$ by (\ref{f19}),
$|d(x,x')-d(y,y')|\leq 2\cdot 2^{-i}$. Therefore,
$|d(x,x')-d(\nu(u_{v_i\langle n,n  \rangle}),\nu(u_{v_i'\langle n,n \rangle}))|
\leq 2\cdot 2^{-i}$.
Since by  (\ref{f7}) $d'(v_i,v_i')=\lim _{n\to\infty}
d(\nu(u_{v_i\langle n,n  \rangle}),\nu(u_{v_i'\langle n,n \rangle}))$,
$|d(x,x')-d'(v_i,v_i')|\leq 2\cdot 2^{-i}$. Then by (\ref{f32}),
\begin{eqnarray}
\label{f29}d(x,x')=d_S((v_0,v_1,\ldots),(v_0',v_1',\ldots))=d_M(\delta_C(q),\delta_C(q'))\,.
\end{eqnarray}
If $\delta(p)=\delta(p')$ then $\delta_C\circ f_N(p)=\delta_C\circ f_N(p')$, hence $f_N$ realizes a single-valued function $f:X\to M$. By (\ref{f29}),
\begin{eqnarray}\label{f34}
d(x,x') & = & d_M(f(x),f(x'))\,,\end{eqnarray}
therefore, $f$ is a $(\delta,\delta_C)$-computable isometric function.
\smallskip

Finally, we show that $f^{-1}$ is $(\delta_C,\delta)$-computable.
Suppose $f(x)=y=\delta_C(q)\in\range(f)$ with $q=(\iota(w_0)\iota(w_1)\ldots)$. Notice that not necessarily $\nu(w_{n+1})\In \nu(w_n)$.
There is some $p'$ such that $x=\delta(p')$ and, by the definition of $f$, $y=f(x)=\delta_C\circ f_N(p')$. Then there are words $w'_n$ such that
$q':=f_N(p')=(\iota(w_0')\iota(w_1')\ldots)$ and $\delta_C(q')=y$.
By (\ref{f17}, \ref{f18}), for these words,
$x\in \nu(w_{n+1}')\In\nu(w_n')$, $ \dm\circ\nu (w_n')<2^{-n}$
and hence $\lim_{n\to\infty}d'(w_n,w_n')=0$.
Since $\delta_C(q)=y=\delta_C(q')$,
\begin{eqnarray}\label{f14}d'(w_n,w_n')&\leq & 2\cdot2^{-n}\ \ \mbox{ for all} \ n\,.
\end{eqnarray}
(Still, for every $w\in\dom(\nu)$ let $((u_{wk},a_{wk}))_{k\in\IN}$ be the sequence satisfying (\ref{f5}, \ref{f6})  for $k=\langle i,n\rangle$ that has been used for defining the pseudometric space ${\bf A}'$.)
Since $x\in\nu(w_n')$ and $\nu(u_{w_n'\langle n,n \rangle})\In \nu(w_n')$,
$d(x,\nu(u_{w_n'\langle n,n \rangle}))\leq \dm(w_n')\leq 2^{-n}$.\\
Suppose $z\in \nu(u_{w_n\langle n,n \rangle})$. Then $d(\nu(u_{w_n\langle n,n \rangle}),z)=0$, hence by (\ref{f22}, \ref{f8},\ \ref{f14})
\[\begin{array}{llll}
d(x,z)&\leq & d(x, \nu(u_{w_n'\langle n,n \rangle})) +
d(\nu(u_{w_n'\langle n,n \rangle}),\nu(u_{w_n\langle n,n \rangle}))\\
&&+ d(\nu(u_{w_n\langle n,n \rangle}),z) + 10\cdot 2^{-n}\ \ \
\\
&\leq& 2^{-n} + d'(w_n',w_n) +10\cdot 2^{-n}\\\
&\leq& 13\cdot 2^{-n}\,.\\
\end{array}\]

Let $((u_i,v_i))_{i\in\IN}$ be the computable numbering of the relation $R$ from
(\ref{f39}) for defining the functions $f_i$ in the proof of Theorem~\ref{t5}.
We prove that $x\in\nu(v)$ iff there are numbers $\;i,n\in\IN$ such that
\begin{eqnarray}\label{f20}
n\geq i+4, \ \ v=v_i,\ \ \mbox{and}\ \ \nu(u_i)\cap \nu(u_{w_n\langle n,n\rangle})\neq\emptyset\,.
\end{eqnarray}
Suppose, $x\in\nu(v)$. There is some $i$ such that $v=v_i$ and $x\in \nu(u_i)\In \nu(v_i)$. There is some $j$ such that $x\in B(x,2^{-j})\In \nu(u_i)$.
Let $n:= \max(j+4,i+4)$. Then for all $z\in \nu(u_{w_n\langle n,n \rangle})$, \
$d(x,z)\leq 13\cdot 2^{-n}<2^{-j}$. We conclude
$\nu(u_{w_n\langle n,n \rangle})\In B(x,2^{-j})
\In \nu(u_i)$,
hence  $\nu(u_i)\cap \nu(u_{w_n\langle n,n\rangle})=\nu(u_{w_n\langle n,n\rangle})\neq\emptyset$ (since $(\forall u)\nu(u)\neq\emptyset$ by assumption).

On the other hand, suppose (\ref{f20}) holds for some $i,n\in\IN$. There is some
$z\in \nu(u_i)\cap \nu(u_{w_n\langle n,n\rangle})$. Then $d(x,z)\leq 13\cdot 2^{-n}$
as shown above. Therefore by $n\geq i+4$,
\[|f_i(x)-f_i(z)|\leq 2^i \cdot d(x,z)\leq 2^i \cdot 13\cdot  2^{-n}<1\,.\]
Since $z\in \nu(u_i)$, $f_i(z)=0$. Therefore, $f_i(x)<1$ hence $x\in\nu(v_i)=\nu(v)$
(see (\ref{f2},\ref{f3},\ref{f4})).

Since $(w,n)\mapsto u_{w,\langle n,n\rangle}$ is computable and $\nu(u)\cap\nu(v)\neq\emptyset$ is r.e., by (\ref{f20}) from $q=\iota(w_0)\iota(w_1)\ldots$ we can compute a list of all $v$ such that $x\in\nu(v)$.  Therefore, $f^{-1}$ is $(\delta_C,\delta)$-computable.
\qq\\

We mention that by (\ref{f34}) the embedding $f$ is an isometric function from the metric space $(X,d)$ into $(M,d_M)$, where $d$ is the metrization of the original $T_3$-space $\bf X$  constructed in the proof of Theorem~\ref{t5}.
Let ${\rm NE}$ abbreviate ``$U\neq \emptyset$ is $\nu$-r.e.''.
The condition in the embedding theorem~\ref{t8} is $({\rm CT_3 +  NE})$. By Theorem~\ref{t7} this implies $\rm SCT_3$.
 We do not know whether $\rm SCT_3$ or $\rm STy$ are sufficient to prove the embedding theorem.

\section{Separation on Product Spaces} \label{secr}
For a computable topological space ${\bf X}=(X,\tau,\beta,\nu)$ and $B\In X$ the subspace
${\bf X}_B=(B,\tau_B,\beta_B,\nu_B)$  of ${\bf X}$ to $B$ is the computable topological space defined
by $\dom(\nu_B):=\dom(\nu)$, $\nu_B(w):=\nu(w)\cap B$ \cite[Section~8]{WG09}.
The separation axioms from Definition~\ref{d4} are invariant under restriction to subspaces.

\begin{thm}\label{t4}
If a computable topological space satisfies  some separation axiom from
Definition~\ref{d4}  then each subspace satisfies this axiom.
\end{thm}

\pproof Straightforward.
\qq\\

The product of two $T_i$-spaces is a $T_i$-space for $i=0,1,2,3$. This is no longer true for some of the computable separation axioms. The product
${\bf X}_1\times {\bf X}_2=\overline{\bf X}=(X_1\times X_2,\overline\tau,\overline\beta,\overline\nu)$
of two computable topological spaces ${\bf X}_1=(X_1,\tau_1,\beta_1,\nu_1)$ and
${\bf X}_2=(X_2,\tau_2,\beta_2,\nu_2)$, defined by  $\overline\nu\langle u_1,u_2\rangle=\nu_1(u_1)\times \nu_2(u_2)$, is again a computable topological space
\cite[Section~8]{WG09}. For the next theorem see Figure~\ref{fig3}.

\begin{thm}\label{t2}
\item\label{t2a} The $SCT_2$-, $WCT_3$-, $CT_3$-, $CTy$- and $SCT_3$-spaces are closed under finite products.
\end{thm}

We consider computability w.r.t. $\nu_i$, $\delta_i$, $\psi^-_i$, $\overline\nu$, $\overline\delta$  and $\overline\psi^-$.

\pproof
Suppose,   ${\bf X}_1$ and ${\bf X}_2$ are $SCT_2$. By \cite[Theorem~7]{Wei10}
, $x_i\neq y_i$ is
$(\delta_i,\delta_i)$-r.e. for $i=1,2$, hence $(x_1,x_2)\neq (y_1,y_2)$ is
$([\delta_1,\delta_2],[\delta_1,\delta_2])$-r.e.,  hence again by Theorem~\ref{t5},
${\bf X}_1\times {\bf X}_2$ is $SCT_2$.

Suppose,   ${\bf X}_1$ and ${\bf X}_2$ are $WCT_3$. Let $(x_1,x_2)\in W_1\times W_2$.
From $x_i$ and $W_i$ we can find $U_i\in\beta_i$ such that $x_i\in U_i\In \overline U_i\In W_i$ (for $i=1,2$). Then $(x_1,x_2)\in U_1\times U_2\In \overline{U_1\times U_2}=
\overline U_1\times \overline U_2\In W_1\times W_2$.

Suppose,   ${\bf X}_1$ and ${\bf X}_2$ are $CT'_3$.  Suppose $(x_1,x_2)\in (W_1,W_2)\in\beta_1\times \beta_2$. From \\
$((x_1,x_2), (W_1,W_2))$ we can compute $x_1$, $x_2$, $W_1$ and $W_2$. Using $t_3'$ for $\bf X_1$ and $\bf X_2$ we can compute $(U_i,B_i)$ such that $U_i\in\beta_i$, $B_i\In X_i$ is closed and $x_i\in U_i\In B_i\In W_i$ ($i=1,2$). Observe that
$(x_1,x_2)\in U_1\times U_2\In B_1\times B_2\In W_1\times W_2$. Form $(U_1,B_1)$ and $(U_2,B_2)$
we can compute $((u_1,u_2),(B_1,B_2))$.

Suppose ${\bf X}_1$ and ${\bf X}_2$ are $CTy'$. We show that ${\bf X}_1\times{\bf X}_2$ is $CTy'$. From $(x_1,x_2)\in W_1\times W_2$ where $W_1\in\beta_1$ and  $W_2\in\beta_2$,  we can compute $x_1$ and $W_1$, where $x_1\in W_1$.
By $\rm CTy'$ for ${\bf X}_1$, from these data we can compute some $U_1\in\beta_1$ and a function $f_1:X_1\to\IR$ such that $x_1\in U_1 \In W_1$ and $f_1$ is $0$ inside $U_1$ and $1$ outside $W_1$. Correspondingly, we can compute $U_2$ and $f_2$ such that
$x_2\in U_2 \In W_2$ and $f_2$ is $0$ inside $U_2$ and $1$ outside $W_2$. From $U_1,U_2$ we can compute $\overline U:=U_1\times U_2$ and $\overline f$ such that $\overline f(y_1,y_2)=\max
(f_1(y_1), f_2(y_2))$. Then $(x_1,x_2)\in U_1\times U_2\In W_1\times W_2$, and
and $\overline f$ is $0$ inside $U_1\times U_2$ and $1$ outside $W_1\times W_2$.

 For $\bf X_i$ ($i=1,2$) let $R_i$ be the r.e. set
and let $r_i$ be the computable function for $SCT_3$ from Definition~\ref{d4}.
By \cite[Lemma~27]{WG09} there is a computable function $h$ such that
$\psi^-_1(p_1)\times \psi^-_2(p_2)=\overline \psi^-\circ h( p_1,p_2)$.
Let
\begin{eqnarray*}
&\overline R  := \{(\langle u_1,u_2\rangle, \langle w_1,w_2\rangle)\mid (u_1,w_1)\in R_1 \wedge
(u_2,w_2)\in R_2\}\,,&\\
&\overline r(\langle u_1,u_2\rangle, \langle w_1,w_2\rangle) :=  h(r_1(u_1,w_1),r_2(u_2,w_2))\,.&
\end{eqnarray*}
A straightforward calculation shows that $\overline R$ is the r.e. set
and  $\overline r$ be the computable function for $SCT_3$ from the definition for the product
${\bf X}_1\times {\bf X}_2$.
\qq\\

The $CT_2$-spaces are not closed under product \cite[Theorem~15]{Wei10}.
Presumably, the $CT_4$-spaces, and hence the $CUr$-spaces, are not closed under product.

\section{Final remarks and thanks}\label{secx}
The list of axioms of computable separation in Definition~\ref{d4} is not exhaustive, there may be other ones. Applications must show which of these axioms are the most natural and useful ones. Many questions about the logical relation between the given axioms have not been answered.

I thank the unknown referees for very careful reading and for giving a number of useful remarks.

\bibliographystyle{plain}
%\bibliography{literature,meinebib}
%\end{document}

\section*{Appendix A (Proof of Theorem~\ref{t6}(\ref{t6c}))}

$\bf SCT_3\Longrightarrow CT_4$:  Let $A,B$ be disjoint closed sets.
Suppose there are sequences of open
sets $V_i,W_i$ and of closed sets $S_i,T_i$ ($i=0,1,\ldots$) such that
\begin{eqnarray}
\label{f43}& V_i\In S_i,\ \ W_i\In T_i\,,\\
\label{f44}& A\In\bigcup_{j\in\IN} W_j,\ \ B\cap T_i=\emptyset\,,\\
\label{f45} &B\In\bigcup_j V_{j\in\IN},\ \ A\cap S_i=\emptyset \,.
\end {eqnarray}
For $i\in\IN$ let
\begin{eqnarray}
\label{f46}  G_i:=W_i\setminus \bigcup_{j\leq i}S_i,&& H_i:=V_i\setminus \bigcup_{j\leq i}T_i\,.
\end {eqnarray}
By (\ref{f44},\ref{f45}),
\begin{eqnarray}
\label{f47}  A\In O_A:=\bigcup _i G_i\,,&&  B\In O_B:=\bigcup _i H_i\,.
\end {eqnarray}
The sets $O_A$ and $O_B$ are open. By (\ref{f46}) for $j\leq i$, \ $G_i\cap S_j =\emptyset$
and so  $G_i\cap V_j =\emptyset$. Therefore, $G_i\cap H_j=\emptyset$ for $j\leq i$.
Similarly, $H_i\cap G_j=\emptyset$ for $j\leq i$. Therefore,
$G_i\cap H_j=\emptyset$ for $i,j\in\IN$ and so $O_A\cap O_B=\emptyset$.

It remains to show that sets $O_A$ and $O_B$ can be computed from $A$ and $B$.
Assume $\psi^-(p)=A$ and $\psi^-(q)=B$. From $p$ and $q$ sequences of pairs $(v_i,v^p_i)$ and $(w_i,v^q_i)$ of words can be computed
such that
\begin{eqnarray*}
\{(u,v)\in R\mid v\ll p\}&=&\{(v_0,v^p_0),(v_1,v^p_1),\ldots \}
\\
\{(u,v)\in R\mid v\ll q\}&=&\{(w_0,v^q_0),(w_1,v^q_1),\ldots \},
\end{eqnarray*}
where  $R$ is the r.e. set from $\rm (SCT_3)$.
For $i\in\IN$ let
\begin{eqnarray*}
V_i:=\nu(v_i),&\ & S_i:=\psi^-\circ r(v_i,v^p_i),\\
W_i:=\nu(w_i),&\ & T_i:=\psi^-\circ r(w_i,v^q_i)\,.
\end{eqnarray*}
Then (\ref{f43},\ref{f44},\ref{f45}) hold true.
By \cite[Theorem~11]{WG09} finite intersection and countable union of open sets can be computed, therefore,  from the $V_i,S_i,W_i$ and $T_i$
the open sets sets $O_A$ and $O_B$ defined in  (\ref{f46},\ref{f47}) can be computed, for which
$A\In O_A$, $B\In O_B$ and $O_A\cap O_B=\emptyset$.
Therefore, the multi-function $t_4$ is $(\psi^-,\psi^-,[\theta,\theta])$-computable.

\section*{Appendic B, Proof of Theorem~\ref{t6}(\ref{t6e})}
$\bf CUr\Longrightarrow CT_4:$ See the proof of Theorem~\ref{t6}.

$\bf CT_4\Longrightarrow CUr:$ We effectivize the classical proof from \cite{Eng89}.
Suppose the space is $CT_4$. Then the multi-function
\begin{eqnarray}
\label{f24}
t:(D,U)\mmto (V,C)\ \ \mbox{for open}\ \ U,V\  \mbox{and closed} \
C,D\ \ \mbox{such that} \ \ D\In V\In C\In U\ \
\end{eqnarray}
is computable. (Find  $(V,W)\in t_4(D,U^c)$ and let $C:=W^c$.)
From closed disjoint sets $A,B$ we compute a family $(V_a,C_a)$, $a\in\IQ\cap[0;1]$, of pairs of sets such that
\begin{eqnarray}
\label{f11}&&V_a\ \ \mbox{is open},\ \ C_a\ \ \mbox{is closed},\ \ V_a\In C_a\,,\\
\label{f12}&&C_a\In V_b\ \ \mbox{if}\ \ a<b\,,\\
\label{f13}&&A\In V_0,\ \ C_1\In B^c\,.
\end{eqnarray}
For this purpose let $i\mapsto r_i$ be a canonical bijective numbering of the rational numbers from the interval $[0;1]$ such that $r_0=0$ and $r_1=1$. Define recursively
\begin{eqnarray}
\label{f25}(V_0,C_0)\in t(A,B^c),\ \  (V_1,C_1)\in t(C_0, B^c),\ \ (V_k,C_k)\in t(C_l, V_m)
\end{eqnarray}
such that $r_l$ is the maximum of the numbers in $\{r_0,\ldots,r_{k-1}\}$ which are less than $r_k$
and $r_m$ is the minimum of the numbers in $\{r_0,\ldots,r_{k-1}\}$ which are greater than $r_k$.
The properties (\ref{f11},\ref{f12},\ref{f13}) can be verified easily.
We define two real valued functions $f_<$ and $f_>$ on $X$ as follows:
\begin{eqnarray}
\label{f26} f_<(x)&:=& \sup (\{a\mid x\not\in C_a\}\cup\{0\})\,,\\
\label{f27} f_>(x)&:=& \inf (\{a\mid x\in V_a\}\cup\{1\})\,.
\end{eqnarray}
If $x\in V_a$ and $b>a$ then $x\in C_b$, hence
$b\leq a$ if $x\in V_a$ and  $x\not\in C_b$.  Therefore, $f_<(x)\leq f_>(x)$.
Suppose, $f_<(x) > f_>(x)$ for some $x$. Then there is some $c\in\IQ$ such that
$f_<(x) >c> f_>(x)$. Moreover, there are some $b>c$ such that $x\not\in C_b$ and
some $a<c$ such that
$x\in V_a$. But  $a<b$ and $x\in V_a$ implies $x\in C_b$. Contradiction.
Therefore, $f:=f_<=f_>$. The function $f$ has value $0$ on $A$ and value $1$ on $B$ and is
continuous \cite{Eng89}.

We show that $f$ can be computed from $A$ and $B$.
Since the function $t$ in (\ref{f24}) is computable, the function
$(A,B)\mmto (V_{r_i},C_{r_i})_i$ is $(\psi,\psi,[\theta,\psi]^\omega)$-computable by (\ref{f25}). By \cite[Theorem~13.2]{WG09},
$x\in U$ for open $U$ is $(\rho,\theta)$-r.e. and $x\not\in C$ for
closed $C$ is $(\rho,\psi)$-r.e. Therefore,  from $(V_{r_i},C_{r_i})_i$ and $x$ by (\ref{f26})
we can list all $a\in \IQ$ such that $a< f(x)$ and by (\ref{f27}) we can list all $a\in \IQ$
such that $a> f(x)$. Therefore, the function $((V_{r_i},C_{r_i})_i,x)\mapsto f(x)$
is $([\theta,\psi]^\omega),\delta,\rho)$-computable, hence
by type conversion  \cite[Lemma~3.3.15]{Wei00} $(V_{r_i},C_{r_i})_i\mapsto f$ is
$([\theta,\psi]^\omega),[\delta\to\rho])$-computable.
Therefore, the space is $CUr$.

\end{document}